\documentclass[11pt]{article}
\title{Vanishing and non-vanishing criteria in Schubert calculus}
\author{Kevin Purbhoo, University of British Columbia
\footnote{Research partially supported by an NSERC scholarship.}}

\usepackage[mathscr]{eucal}
\usepackage{latexsym}
\usepackage{amssymb}
\usepackage{amsthm}
\usepackage{amsmath}
\usepackage{epsfig}

\newcommand{\frn}{\mathfrak{n}}
\newcommand{\frb}{\mathfrak{b}}

\newcommand{\frg}{\mathfrak{g}}

\newcommand{\frsl}{\mathfrak{sl}}

\newcommand{\CC}{\mathbb{C}}
\newcommand{\ZZ}{\mathbb{Z}}
\newcommand{\RR}{\mathbb{R}}

\newcommand{\into}{\hookrightarrow}
\newcommand{\zeroset}{\{0\}}

\newcommand{\strconst}{c_{\pi_1 \ldots \pi_s}}
\newcommand{\triplestrconst}{c_{\pi_1 \pi_2 \pi_3}}
\newcommand{\onetos}{\{1, \ldots, s\}}
\newcommand{\shift}[1]{\boxplus_{#1}}

\newcommand{\tokens}{\mathcal{T}}
\newcommand{\regions}{\mathcal{R}}
\newcommand{\branchproblem}{i^*(\omega_\pi)}

\newcommand{\sqr}[1]{\multicolumn{1}{|c}
{\rule[-8pt]{0pt}{24pt}$#1$}}
\newcommand{\sqrend}[1]{\multicolumn{1}{|c|}
{\rule[-8pt]{0pt}{24pt}$#1$}}

\newcommand{\xothertheoremname}{Theorem}
\newtheorem*{xothertheorem}{\xothertheoremname}

\newtheorem{theorem}{Theorem}
\newtheorem{lemma}{Lemma}[section]
\newtheorem{proposition}[lemma]{Proposition}
\newtheorem{corollary}[lemma]{Corollary}

\newtheorem{example}{Example}[section]
\newtheorem{definition}[example]{Definition}
\newtheorem{remark}[example]{Remark}

\begin{document}
\maketitle

\begin{abstract}
For any complex reductive connected Lie group $G$, many of the structure constants
of the ordinary
cohomology ring $H^*(G/B; \ZZ)$ vanish in the Schubert basis, and the rest 
are strictly 
positive.  We present a combinatorial game, the ``root game'',
 which provides some criteria 
for determining which of the Schubert intersection numbers vanish.  
The definition of the root game is manifestly invariant under automorphisms
of $G$, and under permutations of the classes intersected.
Although these criteria are not proven to cover all cases, in practice 
they work very well, giving a complete answer to the question for 
$G=SL(7,\CC)$.  In a separate paper we show that one of these 
criteria is in fact necessary and sufficient when the classes are pulled
back from a Grassmannian.

More generally If $G' \into G$ is an inclusion of complex reductive 
connected Lie 
groups, there is an induced map $H^*(G/B) \to H^*(G'/B')$ on the cohomology of 
the homogeneous spaces.  The image of a Schubert class under this
map is a positive sum of Schubert classes on $G'/B'$.  We investigate 
the problem of determining which Schubert classes appear with non-zero
coefficient.  This
is the vanishing problem for branching Schubert calculus, which 
plays an important role in representation theory and symplectic 
geometry, as shown in [Berenstein-Sjamaar 2000].  The root game
generalises to give a vanishing criterion and a non-vanishing
criterion for this problem. 
\end{abstract}


\section{Introduction}
\label{sec:intro}

In this paper we introduce some techniques for studying vanishing
problems in Schubert calculus.  The most basic and famous such problem 
concerns
the cohomology ring of a generalised flag manifold $G/B$: we would
like to determine combinatorially which of the structure
constants for $H^*(G/B)$ are non-zero.  We refer to this as the
vanishing problem for multiplication in Schubert calculus.  However, 
the techniques we introduce here 
apply in a more general context, namely to the vanishing problem for
{\em branching Schubert calculus}---discussed below---of which the
multiplication problem is a special case.
Although it is possible to calculate any structure constants for
these problems explicitly (e.g. using Schubert polynomials \cite{BH}),
the known methods involve alternating sums, and thus
provide little insight into the question of which terms vanish.
A complete combinatorial solution to either of these problems is still
not known.

Our first objective is to provide some vanishing and
non-vanishing criteria for intersection numbers of Schubert varieties on 
$G/B$.  
Geometrically, the problem is this: given $s \geq 3$ Schubert varieties in 
general position, determine whether or not their intersection is empty.
If we know the Schubert intersection numbers we also implicitly have
the Schubert structure constants for $H^*(G/B)$ (from the Poincar\'{e} 
pairing), thus this also addresses the vanishing problem for 
multiplication.

In Section \ref{sec:igame}, we introduce the {\em root game} which can 
often give information about a Schubert intersection number.  
In some circumstances the root game will tell us that the intersection 
number is $0$
(Theorem \ref{v-thm:lose}); in other circumstances, the game will tell us
that the intersection number is at least $1$ (Theorem \ref{v-thm:win}).
Unfortunately, in a few cases, the root game gives no information; 
remarkably though, for $G=SL(n)$, $n \leq 7$ we have confirmed by computer
that all of these remaining cases have intersection 
number $0$.  

The rules of the root game are manifestly symmetric under
permutations of the classes intersected, as well as under automorphisms
of $G$.  Furthermore, once the game has been fully internalized, it
is highly amenable to computations by hand.
 
Our second objective is to show that the main results
hold in an
even more general setting, which we call {\em branching Schubert calculus}.
Let $i:G' \into G$ be an inclusion of complex reductive connected
Lie groups.  Choose Borel subgroups $B' \subset G$ and $B \subset G$ 
such that $i(B') \subset B$.
Then we obtain
an inclusion $i: G'/B' \into G/B$ (which we also denote by $i$, in a
mild abuse of notation).  Hence there is a map
on cohomology $i^*:H^*(G/B) \to H^*(G'/B')$.  The problem of
branching Schubert calculus is to determine the map $i^*$ in the
Schubert basis, i.e.  given a Schubert class $\omega \in H^*(G/B)$
we would like to express $i^*(\omega) \in H^*(G'/B')$ in the Schubert
basis of the latter.

The coefficients which appear in such an expression are always non-negative
integers.  Although there are formulae for these integers, it is not 
known how to determine them combinatorially, or even how to
determine which terms appear.  In Section \ref{sec:bgame}, 
we investigate the latter problem, and obtain some widely applicable
criteria for determining which terms appear.

The vanishing problem for branching Schubert calculus 
generalises the vanishing problem for multiplication
in Schubert calculus: if $i: G' \into G = G' \times G'$
is the diagonal inclusion, then the map $i^*$ is just the cup product
in cohomology.  Similarly, multiplication
of more than two terms comes from considering the diagonal inclusion
$G' \into G' \times \cdots \times G'$.

Our motivation for this work comes from \cite{BS}, in which Berenstein 
and Sjamaar use the vanishing problem for branching
Schubert calculus to answer questions in symplectic geometry and
representation theory.  Let $K'$ and $K$ be the maximal compact subgroups
of $G'$ and $G$ respectively.   Berenstein and Sjamaar use the vanishing 
problem for branching Schubert calculus to calculate the $K'$ moment
polytope of a $K$-coadjoint orbit.  They show that each non-vanishing
branching coefficient gives rise to an inequality satisfied by the
moment polytope.  Moreover, all together,
the complete list of non-vanishing branching coefficients gives
a sufficient set of inequalities for this polytope.  

This symplectic problem 
is to be equivalent to an asymptotic version of a fundamental 
representation theory
question, as shown in \cite{H, GS} 
(for more of this picture see also 
\cite{GLS}).  
Let $\lambda$ and $\mu$ be dominant weights for $G$
and $G'$ respectively.  Let $V_\lambda$ denote the irreducible 
$G$-representation with highest weight $\lambda$; similarly let
$V'_\mu$ denote the irreducible $G'$-representation with highest 
weight $\mu$.  When $V_\lambda$ is decomposed as a $G'$-module,
it is a basic question whether a component of type $V'_\mu$ appears.
The asymptotic version of this problem is the following: does there 
exists a positive integer $N$, such that the $G$-module 
$V_{N\lambda}$ has a component of type $V'_{N\mu}$, when decomposed
as a $G'$-module?  The answer is yes if and only if
the point $\mu$ lies in the $K'$-moment polytope for the $K$-coadjoint
orbit through $\lambda$.  Thus the non-vanishing branching coefficients
give an answer to this asymptotic representation theory question as well.

In studying the vanishing problem for branching Schubert calculus,
we will actually be considering the following apparently
simpler problem: 
determine which Schubert classes are in the kernel of $i^*$.  
While this may at first seem to be a vast simplification, it is
in fact equivalent to the original problem, as shown by Proposition
\ref{prop:kernel}.
In the case of vanishing for multiplication
of Schubert classes, this is a familiar fact: we
can determine which structure constants of the cohomology ring are
zero, based on the which triple products vanish.

The paper begins with a discussion of the geometry underlying the root games
(section \ref{sec:geom}).  The basic idea is to use Kleiman's Bertini
theorem \cite{Kl} to reduce the vanishing problem to a transversality
problem in the tangent space to a point in $G/B$.  Given an
intersection in the tangent space, we attempt to show it is transverse
by degenerating it to a position where transversality is easily verifiable.
If this is possible, we can conclude that
that a corresponding Schubert class is not in the kernel of $i^*$.

The degenerations in question can be encoded combinatorially; doing so
gives the root game.
In Section \ref{sec:igame}, we introduce the root game for Schubert
intersection numbers.  Section \ref{sec:bgame} contains the more
general root game for branching Schubert calculus, and proofs of the
main theorems.  Ultimately it is the proof of Theorems \ref{b-thm:lose}
and \ref{b-thm:win}, which tie the combinatorics into the geometry.

We refer the reader to \cite{F} for a general reference on type $A$
Schubert calculus, and to \cite{BH} for the other classical Lie groups.

The author is deeply grateful to Allen Knutson, for providing
a lot of helpful feedback on this paper.


\section{Geometry of vanishing problem for branching Schubert calculus}
\label{sec:geom}


\subsection{Conventions}
\label{sec:conventions}
Given $i:G' \into G$, an inclusion of complex reductive connected
Lie groups, we wish to study the map $i^*:H^*(G/B) \to H^*(G'/B')$.
First, we need to demonstrate that the derived map $i:G'/B' \into G/B$ 
always exists.
\begin{proposition}
Given $i:G' \into G$ there exist Borel subgroups
$B' \subset G'$ and $B \subset G$ such that $i(B') \subset B$.
\end{proposition}

\begin{proof}
Choose a Borel subgroup $B_0 \subset G$, and  consider the 
$G'$-orbits on $G/B_0$ of minimal dimension.  Each such orbit is
closed, therefore, compact, and so is $G'/P$ for some parabolic subgroup
$P \subset G'$.  
Choose a point ${x_0}$ on such an orbit.  The stabiliser of ${x_0}$ inside
$G'$, $G'_{x_0}$, is conjugate to $P$, whereas the stabiliser of ${x_0}$ 
inside $G$, $G_{x_0}$, is conjugate to $B_0$.  Thus $G'_{x_0} \subset G_{x_0}$ is solvable, 
but $G'/G'_{x_0}$ is compact, hence $G'_{x_0}$ is a Borel subgroup of
$G'$.  We take $B = G_{x_0}$ and $B' = G'_{x_0}$.  
\end{proof}

Let $T' \subset B'$ be a maximal torus of $G'$.  Extend its image $i(T')$
to a maximal torus $T \subset B$ of $G$.
Let $N'$ and $N$ denote
the corresponding unipotent subgroups of $B'$ and $B$.  Of 
course, $i(N') \subset N$.  Henceforth we will simply view $T'$, $N'$, $B'$
as subgroups of $T$, $N$, $B$ respectively.

Let $\Delta$ denote the root system of $G$, and $\Delta'$ 
the root system of $G'$.  The positive and negative roots of
$\Delta$ (with respect to the choice of $B$) are denoted 
$\Delta_+$ and $\Delta_-$ respectively.
For each root $\alpha \in \Delta$, we fix a basis vector $e_{\alpha}$
for the corresponding root space in $\frg$.
Likewise, for each root $\beta \in \Delta'$, we fix a basis vector 
$e'_{\beta}$ for the corresponding root space in $\frg'$.

The tangent spaces to ${x_0}$ in $G/B$ and $G'/B'$ are naturally
are $\frg/\frb$ and $\frg'/\frb'$ respectively.  Thus linearising
$i$ gives a a natural inclusion of tangent spaces
$\frg'/\frb' \into \frg/\frb$.
We use the Killing form to identify $\frn$ with 
$(\frg/\frb)^*$.
Similarly, we identify the dual of $\frg'/\frb'$ with $\frn'$.
Thus we obtain a
linear map
$$ \phi : \frn \to \frn'$$
which is adjoint to the inclusion of tangent spaces
$\frg'/\frb' \into \frg/\frb$. 
Essentially $\phi$ encodes all the information about the inclusion 
$G' \into G$.

Note that since ${x_0}$ is a $T'$-fixed point, the map $\phi$ is 
$T'$-equivariant.  Thus, it takes the $T$-weight spaces to 
$T'$-weight spaces, and induces a map
$$\hat \phi : \Delta_+ \to \Delta'_+ \cup \zeroset$$
defined by the rule
$$
\hat \phi(\alpha) = 
\begin{cases}
 0, &\text{if $\phi(e_\alpha) = 0$} \\
 \beta, &\text{where $0 \neq \phi(e_\alpha)$ is in the $\beta$-weight space.}
\end{cases}
$$
In Section \ref{sec:bgame} we will need to consider
subsets $\tokens \subset \Delta_+$ with the following properties.

\begin{definition}
\label{b-def:injective}
Suppose $\tokens \subset \Delta_+$ satisfies
\begin{verse}
1. $0 \notin \phi(\tokens)$, and \\
2. $\hat \phi|_\tokens$ is injective.
\end{verse}
We call such a subset $\tokens$ {\bf injective}.  Equivalently $\tokens \subset
\Delta_+$ is injective if $\phi|_{\langle e_\alpha\,|\,\alpha \in \tokens\rangle}$
is an injective linear map.
\end{definition}


\subsection{Schubert varieties}

Let $W = N(T)/T$ be the Weyl group of $G$.  
For $\pi \in W$, let $[\pi]$
denote the corresponding $T$-fixed point on $G/B$, and let $\tilde{\pi}$
denote some lifting of $\pi \in W$ to an element of 
$N(T) \subset G$.

Let $w_0$ denote the long element in $W$.  For $\pi \in W$, let $\pi' =
w_0 \pi$.
To each $\pi \in W$ we associate the Schubert cell
$X^\circ_{\pi}  = B \cdot [\pi']$, the $B$-orbit through
$[\pi]$ in $G/B$.  Its closure 
$X_{\pi}  = \overline{B \cdot [\pi']}$,
is the Schubert variety.  
(This definition is slightly non-standard: it is more common
to define $X_{\pi} = \overline{B_- \cdot [\pi]}$, where
$B_-$ is an opposite Borel.  Our $X_\pi$ is a translation of the more
standard one by $\tilde w_0$.)
According to these conventions $X_1 = G/B$ 
(where $1 \in W$ represents the identity element) 
and $X_{w_0} = \{x_0\}$.
In general $X_\pi$ is a complex subvariety of $G/B$ whose codimension
\footnote{All dimensions/codimensions are over $\CC$, unless otherwise 
specified.} 
is length of $\pi \in W$ (denoted $\ell(\pi)$).

Let $\omega_\pi$ denote the cohomology class Poincar\'{e} dual
to the homology class of the
Schubert variety $X_\pi$, that is the class such that
$$\int_{X_\pi} \sigma = \int_{G/B} \omega_\pi \cdot \sigma$$  
for all $\sigma \in H^*(G/B)$.  Since the the codimension 
of $X_\pi$ is $\ell(\pi)$, $\omega_\pi$ is a cohomology class of degree
$2\ell(\pi)$.  The class $\omega_1 \in H^*(G/B)$ is the (multiplicative)
identity element.

The following proposition shows that the vanishing
problem for branching Schubert calculus is equivalent to the problem
of determining whether $\omega_\pi \in \ker i^*$.

\begin{proposition}
\label{prop:kernel}
Given $j:G' \into G''$ an inclusion of complex reductive connected
groups, let $G = G' \times G''$ and
$i = (id \times j) \circ \delta:G' \into  G$, where 
$\delta:G' \to G' \times G'$ is the diagonal map.
Let $i^*: H^*(G/B) \cong
H^*(G'/B') \times H^*(G''/B'') \to H^*(G'/B')$ be the induced map
on cohomology.
A Schubert class $\sigma \in H^*(G'/B')$
appears in the expansion of $j^*(\omega)$ if and only if 
$(\sigma^\vee, \omega) \notin \ker i^*$, and 
$\deg \sigma^\vee + \deg \omega = \dim_\RR G'/B'$.  
Here  $\sigma^\vee$ is
the Schubert class dual to $\sigma$ under the Poincar\'{e}
pairing .
\end{proposition}

\begin{proof}
Consider the integral
$$\int_{G'/B'} \sigma^\vee \cdot j^*(\omega).$$
If this integral is non-zero, then 
$\sigma$ appears in the
expansion of $j^*(\omega)$ (with coefficient equal to 
$\int_{G'/B'} \sigma^\vee \cdot j^*(\omega)$);
otherwise it does not.

Since $i=(id \times j)\circ \delta$, we have
$i^*(\sigma^\vee, \omega) =
\delta^*(\sigma^\vee, j^*(\omega)) =
\sigma^\vee \cdot j^*(\omega)$,
thus 
$$
\int_{G'/B'} \sigma^\vee \cdot j^*(\omega) =
\int_{G'/B'} i^*(\sigma^\vee, \omega).$$ 
The second integral is clearly non-zero
if and only if  $(\sigma^\vee, \omega) \notin \ker i^*$ and 
$\deg \sigma^\vee + \deg \omega = \dim_\RR G'/B'$.  
\end{proof}

Thus, to solve the vanishing problem for branching Schubert calculus
for $j:G' \into G''$,
it is sufficient to know whether $i^*(\sigma^\vee, \omega) = 0$,
for any given $(\sigma^\vee, \omega) \in H^*(G/B)$.

Henceforth we shall be investigating
the question of whether $\branchproblem = 0$, for $\pi \in W$.
We will assume that $\pi \in W$ is an element whose length 
$\ell(\pi) \leq \dim G'/B'$:
if $\ell(\pi) > \dim G'/B'$ then $\branchproblem = 0$ for dimensional 
reasons.
We are primarily interested in the case where $\ell(\pi) = \dim G'/B'$,
however except where specified otherwise, 
everything in this paper holds for all $\pi \in W$.


\subsection{The multiplication problem}
A special and particularly important case is the vanishing problem
for multiplication of Schubert calculus.  As mentioned, in the
introduction, this corresponds to the diagonal inclusion
$G' \into G = G' \times \cdots \times G'$ ($s$-factors).  
 
In this case, a Schubert class 
$\omega_\pi \in H^*(G/B)$ can be regarded as an $s$-tuple
of Schubert classes 
$(\omega_{\pi_1}, \ldots, \omega_{\pi_s}) \in (H^*(G'/B'))^s$. The
map $i^*: H^*(G/B) \to H^*(G'/B')$ gives the product of
these Schubert classes in $H^*(G'/B')$:
$$\branchproblem = i^*(\omega_{\pi_1}, \ldots, \omega_{\pi_s})
=\omega_{\pi_1} \cdots  \omega_{\pi_s}.$$
Thus the problem of determining when $\branchproblem \neq 0$ becomes
the question of which collections of Schubert
classes on $G'/B'$ have non-vanishing product.

We are most interested in the case where 
$\ell(\pi) = \sum \ell(\pi_i) = \dim(G'/B')$.  In this case we are investigating
the Schubert intersection numbers
$\strconst$ defined by
$$\strconst =
\int_{G'/B'} \omega_{\pi_1} \cdots \omega_{\pi_s}.$$
The triple Schubert intersection numbers
$\triplestrconst$ are particularly important, as they are the 
Schubert structure constants of the cohomology
ring $H^*(G/B)$.   Indeed, if we write
$$\omega_{\pi_1} \cdot \omega_{\pi_2} = \sum_{\rho \in W}
c_{\pi_1 \pi_2}^{\rho} \omega_{\rho}$$
then
\begin{align*}
\triplestrconst &=
\int_{G/B} \omega_{\pi_1} \cdot \omega_{\pi_2} \cdot \omega_{\pi_3} \\
&= \int_{G/B} \sum_{\rho \in W} c_{\pi_1 \pi_2}^{\rho}
\omega_{\rho} \cdot \omega_{\pi_3} \\
&= c_{\pi_1 \pi_2}^{w_0 \pi_3}.
\end{align*}


\subsection{Tangent space methods}

The main idea behind the results in this paper is to use Kleiman's
theorem \cite{Kl} to translated problems of intersection theory on $G/B$ 
into transversality
questions on the tangent space to $G/B$.  Tangent space methods have
been used elsewhere in the literature, perhaps most notably in
Belkale's geometric proof of the Horn conjecture \cite{B}.  Our main
lemma (Lemma \ref{b-lem:linalg}) generalises some of these ideas.

\begin{lemma} \label{lem:tangentspaces}
Let $x \in G/B$.  The following conditions are equivalent:
\begin{enumerate}
\item $\branchproblem \neq 0$,
\item There exist $g_1,g_2 \in G$ such that 
$x \in g_1 X^\circ_\pi \cap g_2 G'/B'$, and the tangent spaces
$T_x g_1 X_\pi$ and $T_x g_2 G'/B'$ are transverse linear subspaces of
$T_x G/B$.
\end{enumerate}
\end{lemma}

\begin{proof}
We apply Kleiman's Theorem to the $G$-homogeneous space $G/B$ and
its subvarieties $X_\pi$ and $G'/B'$.
Consider the intersections
$$I_{g_1,g_2} = g_1 X_\pi \cap g_2 G'/B' $$
and 
$$I^\circ_{g_1,g_2} = g_1 X^\circ_\pi \cap g_2 G'/B'\,.$$
If $g_1, g_2$ are generic elements of $G$, Kleiman's theorem tells
us that a generic point $\tilde x$ of $I_{g_1,g_2}$ is a smooth point of 
$g_1 X_\pi$ which can be assumed to lie in $g_1X^\circ_\pi$; 
moreover the varieties 
$g_1 X^\circ_\pi$ and
$g_2 G'/B'$ are transverse at $\tilde x$.  In particular $I_{g_1,g_2}$ 
is generically
reduced and equidimensional.  If $I_{g_1,g_2}$ is zero-dimensional then
$I_{g_1,g_2}$ is finite, with cardinality 
$\#(I_{g_1,g_2}) = \int_{G'/B'} \branchproblem$.
More generally $I_{g_1,g_2}$ defines a homology class in $G'/B'$ which is
Poincar\'{e} dual to to the cohomology class $\branchproblem$.
In particular, we have that $\branchproblem \neq 0$ if and only
if 
$I_{g_1,g_2}$ (or equivalently $I^\circ_{g_1,g_2}$) 
is nonempty for generic $(g_1,g_2) \in G \times G$.

Let $A = \{(g_1,g_2)\ |\ I^\circ_{g_1,g_2} \neq \emptyset \}$.
(Note $G' \times G' \subset A$ so $A$ is always non-empty.)
We have just shown 
$\bar A = G \times G$ if and only if $\branchproblem \neq 0$.
Let $(g_1,g_2)$ be a generic point of $A$ (If $A$, is reducible choose
any component), and $\tilde x \in I_{g_1,g_2}$ be
a generic point of $I_{g_1,g_2}$.  If $\branchproblem \neq 0$ then
$(g_1,g_2)$ is in fact a generic point of $G \times G$ and so
the varieties $g_1 X_\pi$ and $g_2 G'/B'$ are transverse at $\tilde x$.
However, note that the set 
$$\{(g_1,g_2)\ |\ \text{$g_1 X^\circ_\pi$ 
and $g_2 G'/B'$ have a transverse point of intersection}\}$$ 
is necessarily open (intuitively this is because a transverse intersection
remains transverse under perturbation).  Thus, conversely, if
$g_1 X_\pi$ and $g_2 G'/B'$ are transverse at $\tilde x$, then 
$\bar A = G \times G$ and hence $\branchproblem \neq 0$.

Finally, since $G$ acts transitively on $G/B$, we can find $g \in G$ such
that $g\tilde x =x$.  Then 
$g_1 X_\pi$ and $g_2 G'/B'$ are transverse at $\tilde x$
iff $gg_1 X_\pi$ and $gg_2 G'/B'$ are transverse at $x$.
This completes the proof.
\end{proof}

Lemma \ref{lem:tangentspaces} is still not concrete enough
for our purposes.  We reformulate it as follows.

For $a \in N$, let $a \cdot \ :\frn \to \frn$ denote the adjoint action
of $N$ on its Lie algebra.
Let $Q \subset \frn$ be
the subspace generated by the $e_\alpha$ such that
$\alpha \in \Delta_+$ and  $\pi^{-1} \cdot \alpha \in \Delta_-$.
Equivalently,
$$Q = \frn \cap (\pi \cdot \frb_-).$$

\begin{lemma}
\label{b-lem:linalg}
The following are equivalent:
\begin{enumerate}
\item $\branchproblem \neq 0$.
\item $\phi|_{a \cdot Q}$ is injective for some $a \in N$.
\item $\phi|_{a \cdot Q}$ is injective for generic $a \in N$.
\end{enumerate}
\end{lemma}

The tangent space to $G/B$ at $x_0$ is naturally $\frg/\frb$.
We identify the cotangent space $(\frg/\frb)^*$ with $\frn$ 
using the Killing form.
Under these identifications, $Q^\perp \simeq ((\pi' \cdot \frb) + \frb)/\frb$.
The subspace $a \cdot Q \subset \frn$ is identified with the conormal 
space at the point $x_0$ to a translated Schubert variety $g \cdot X_\pi
\ni x_0$.  Thus Lemma \ref{b-lem:linalg} is 
essentially a dual statement to Lemma \ref{lem:tangentspaces}.

\begin{proof}
The equivalence of conditions 2 and 3 is clear, as the maps
$\phi|_A$ are
injective for a Zariski open set of subspaces $A$.  

To show
the equivalence of 1 and 3, 
we use Lemma \ref{lem:tangentspaces} with the point $x=x_0$.

We have $x_0 \in g \cdot X^\circ_\pi$ if and only if 
$g =b_1 (\tilde \pi')^{-1}$ for some $b_1 \in B$, and
$x_0 \in g' G'/B'$ if and only if $g' = b_2 h$, for $b_2 \in B$,
$h \in G'$.  Put $b = b_2^{-1}b_1$, and write $b=at$, with
$a\in N$, and $t \in T$.

Then,
\begin{align*}
T_{x_0} g X_\pi \cap T_{x_0} g' G'/B' 
&= b_2 \cdot \big(
b \cdot T_{x_0} (\tilde\pi')^{-1} X_{\pi} \cap 
T_{x_0} G'/B' \big) \\
&= b_2 \cdot \big(
(b \cdot (\pi' \cdot \frb) + \frb)/\frb \cap \frg'/\frb' \big) \\
&= b_2 \cdot \big(
(b \cdot Q)^\perp \cap \frg'/\frb \big) \\
&= b_2 \cdot \big(
(a \cdot Q)^\perp \cap \frg'/\frb \big)
\end{align*}

The transversality of the intersection 
$(a \cdot Q)^\perp \cap \frg'/\frb$
is precisely the dual statement to condition (2).
\end{proof}

Applied to the multiplication problem, Lemma \ref{b-lem:linalg}
reduces to the following.

\begin{corollary}
\label{cor:v-linalg}
Let $Q_i = \frn' \cap (\pi_i \cdot \frb'_-)$ be the subspace of $\frn'$
whose weights are the inversion set of $\pi_i$.  Then the following are
equivalent:
\begin{enumerate}
\item $\strconst \neq 0$,
\item 
The sum of subspaces $a_1 \cdot Q_1 + \cdots + a_s \cdot Q_s$ is 
a direct sum, for generic choices of $a_i \in N$.
\end{enumerate}
\end{corollary}


\subsection{Necessary conditions for vanishing}
Our first consequence of Lemma \ref{b-lem:linalg} is the vanishing
criterion.

\begin{lemma}\label{lem:geomlose}
Let $S \subset \frn$ be an $N$-submodule of $\frn$.  If 
$\dim \phi(S) < \dim(Q \cap S)$, then $\branchproblem =0$.
\end{lemma}

\begin{proof}
As $S$ is $N$-invariant, we have that
$$\dim((a \cdot Q) \cap S) = \dim (Q \cap S) > \dim \phi(S)$$
for all $a \in N$.
It follows that $\phi|_{((a \cdot Q) \cap S)}$ is not injective, and
thus $\phi|_{a \cdot Q}$ is not injective.  Therefore,
by Lemma \ref{b-lem:linalg}, $\branchproblem = 0$.
\end{proof}

Moreover, if we take $S$ to be a $B$-submodule of $\frn$ then there are
only finitely many possibilities, and we can readily
calculate the dimensions of $\phi(S)$ and $Q \cap S$ combinatorially.
This is essentially the content of Theorem \ref{b-thm:lose}.

\begin{remark}
From Lemma \ref{lem:geomlose}, it is possible to rederive the
necessary Horn inequalities for non-vanishing of Schubert calculus
on Grassmannians.  For more on this picture, see \cite{P-horn}.
\end{remark}


\subsection{Degenerating $Q$}
\label{sec:degenerate}

To show $\branchproblem \neq 0$, by Lemma \ref{b-lem:linalg},
it is enough to exhibit a subspace $U= a \cdot Q$ in the $N$-orbit 
through the subspace $Q$ such that $\phi|_U$ is injective.  
Actually, because the set 
$$\{U \in Gr(\dim G'/B', \frn)\ |\ \phi|_U \text{is injective}\}$$
is open, we can take $U$ to be in the closure of
the $N$-orbit through $Q$.  Note that since $Q$ is a $T$-fixed subspace of 
$\frn$, the $B$-orbit through $Q$ coincides with the $N$-orbit through $Q$.  

The idea behind obtaining sufficient conditions is to look for
a $T$-fixed subspace of $\frn$, $U \in \overline{B \cdot Q}$, such that
$\phi|_U$ is injective.  We can think of the search for a suitable $U$
as a process.  Beginning with the $T$-fixed subspace
$Q \subset V$ we degenerate to another $T$-fixed subspace $U
\in \overline{B \cdot Q}$. If $\phi|_U$ is not injective, we can 
degenerate further inside $\overline{B \cdot U}$, until a suitable subspace
is found.

Let $V=\frn$ or any $B$-module subquotient of $\frn$.
Let $V'=\frn$ or any $B'$-module subquotient of $\frn'$.
Suppose we have a $B'$-equivariant map $\psi: V \to V'$

Let $Gr(V)$ denote the disjoint union of all Grassmannians 
$$Gr(V) = \coprod_{l=0}^{\dim V} Gr_l(V).$$
Since $V$ has a $B$-action, so does $Gr(V)$.

Let $U \in Gr(V)$ be a subspace of $V$.
We call the quadruple $(U,V,V', \psi)$ {\bf good} if there is a point
$\tilde U \in \overline{B \cdot U}$ such that
$\psi|_{\tilde U}:\tilde U \to V'$ is an injective linear map.
Note that the set of ${\tilde U} \in Gr(V)$ with $\phi|_{\tilde U}$ 
injective is Zariski open in $Gr(V)$.  Thus, equivalently, 
$(U,V,V', \psi)$ is good if there exists
$\tilde U \in B \cdot U$ such that
$\psi|_{\tilde U}:{\tilde U} \to V'$ is an injective linear map.

In the language of good quadruples, Lemma \ref{b-lem:linalg} states that
$\branchproblem = 0$ if and only if the quadruple $(Q,\frn, \frn', \phi)$
is good.

\subsubsection{Moving between fixed points}
For any $T$-representation $U$ with distinct weights, let
$\Gamma(U)$ denote the set of weights of $U$.

To every $\beta \in \Delta_+$, we can associate a one dimensional
unipotent subalgebra $N_\beta \subset N$, whose Lie algebra 
$\frn_\beta$ is $T$-invariant with weight $\beta$.  $N_\beta$
is isomorphic to the additive Lie group $\CC$.  Let 
$\theta_\beta:\CC = N_\beta \into N$ denote the inclusion of
groups, $\theta_\beta(t) = \exp te_\beta$.

The following proposition is a triviality, yet it is at the very
heart of the root game.

\begin{proposition}\label{prop:geommove}
Let $U \in Gr(V)$, and let 
$U^1 = \lim_{t  \to \infty} \theta_\beta(t) \cdot U$.
If $(U^1,V,V',\psi)$ is good, then  $(U,V,V',\psi)$ is good.
\end{proposition}

\begin{proof}
The point $(U^1,V,V',\psi)$ lies in the closure of 
$\overline{B \cdot U}$.  Thus if there exists 
$\tilde U \in \overline{B \cdot U^1}$ such that $\psi|_{\tilde U}$
is injective, then $\tilde U$ also lies in $\overline{B \cdot U}$.
\end{proof}

\begin{remark}
In particular if $\psi|_{U^1}$ happens to be injective then
$(U,V,V',\psi)$ is good.  Otherwise we can attempt to apply
Proposition \ref{prop:geommove} recursively to 
$(U^1, V, V', \psi)$, 
to show that $(U^1,V,V',\psi)$ is good and hence that
$(U,V,V',\psi)$ is good.
\end{remark}
 
Suppose now that $U$ is a $T$-fixed
point of $Gr(V)$.
We show that $U^1$ is a $T$-fixed point of $Gr(V)$ and 
calculate the weights $\Gamma(U^1)$ in terms 
of $\Gamma(U)$.

\begin{definition}\label{def:shift}
Call an element $\alpha \subset \Gamma(U)$ {\bf $\beta$-shiftable},
if there is a positive integer $k$ such that $\alpha + k\beta
\in \Gamma(V) \setminus \Gamma(U)$.  
Let $\Gamma(U) \shift{\Gamma(V)} \beta$
denote the set
$$\{\alpha+\beta\ |\ \text{$\alpha$ is $\beta$-shiftable}\}
\ \cup\ \{\alpha\ |\ \text{$\alpha$ is not $\beta$-shiftable}\}
$$
\end{definition}

\begin{lemma}\label{lem:rootsmove}
Let $U^1 = \lim_{t  \to \infty} \theta_\beta(t) \cdot U$.  Then
$U^1$ is a $T$-fixed point of $Gr(V)$ and
$\Gamma(U^1)= \Gamma(U) \shift{\Gamma(V)} \beta$.
\end{lemma}

\begin{proof}
Let $\bar e_\alpha \in V$ be a vector with weight $\alpha$.
Since the weights of $V$ are distinct, we can represent $U$ as
$[\bar{e}_{\alpha_1} \wedge \ldots \wedge \bar{e}_{\alpha_l}]$,
and $U^1$ as
$[\bar{e}_{\alpha'_1} \wedge \ldots \wedge \bar{e}_{\alpha'_l}]$,
via the Pl\"{u}cker embedding $Gr(V) \into P(\bigwedge^*V)$.
Now
\begin{align*}
\theta_{\beta}(t) \cdot U
&= \theta_{\beta}(t) \cdot
         [\bar e_{\alpha_1} \wedge \ldots \wedge \bar e_{\alpha_l}],
         \\
&= \big[(\bar e_{\alpha_1} + t(e_\beta \cdot \bar e_{\alpha_1}))
             \wedge \ldots \wedge
   (\bar e_{\alpha_l} + t(e_\beta \cdot \bar e_{\alpha_l})) \big] \\
&= 
\sum_{C \subset \{1, \ldots, l\}}
t^{|C|} \pm \bigwedge_{i \in C} e_\beta \cdot \bar e_{\alpha_i}
\wedge \bigwedge_{i \in C^c} \bar e_{\alpha_i}
\end{align*}
(here $e_\beta\cdot$ is the action of $\frn_\beta$ on 
$V$ induced from the adjoint action).
Now up to a non-zero constant multiple,
$$e_\beta \cdot \bar e_{\alpha_i} = 
\begin{cases}
\bar e_{\alpha_i + \beta}, &\text{if $\alpha_i + \beta \in \Gamma(V)$} \\
$0$, &\text{otherwise.}
\end{cases}
$$ 
This is a property of the adjoint representation which $V$, as subquotient
of the adjoint representation, inherits.

We see that a summand is non-zero only if $\{\alpha_i\ |\ i \in C\}$ 
is a subset of the set of $\beta$-shiftable weights of $\Gamma(U)$.  
In the limit as $t \to \infty$, the only term which
survives is the one with the highest power of $t$, which is precisely
$$[\pm t^{\text{\#$\beta$-shiftable weights}}
      \bar e_{\alpha'_1} \wedge \ldots \wedge \bar e_{\alpha'_l}].$$
\end{proof}

\subsubsection{Splitting into two smaller problems}

Let $S \subset V$ be an $B$-submodule, and $S' \subset V$ an
$B'$-submodule.  Suppose that $\psi(S) \subset S'$.
Let $q:V \to V/S$ and $q':V' \to V'/S'$ denote the quotient maps.  

From the quadruple $(U, V, V', \psi)$ and the submodules $S,S'$, we 
obtain two induced quadruples: they are
$(U \cap S, S, S', \psi|_S)$, and 
$(q(U), V/S, V'/S', \psi_q)$, where $\psi_q:=q'\circ \psi \circ q^{-1}$.
(Note that $\psi_q:V/S \to V'/S'$ is well defined.)

\begin{proposition}\label{prop:geomsplit}
If 
$(U \cap S, S, S', \psi|_S)$, and 
$(q(U), V/S, V'/S', \psi_q)$
are both good, then $(U, V, V', \phi)$ is good.
\end{proposition}

\begin{proof}
Let $p:Gr(V) \to Gr(S)$ be the map $p(U) = U \cap S$, and note
that $q$ also defines a similar map $q: Gr(V) \to Gr(V/S)$.
Note that $p$ and $q$ are not continuous everywhere, but since
$S$ is a $B$-submodule, they are $B$-equivariant and continuous
on $B$-orbits.

Define
$$g(U,V,V',\psi) 
:= \{\tilde U \in B \cdot U \subset Gr(V)
\ |\ \text{$\psi|_{\tilde U}$ is injective}\}.
$$
Let $g_p = g(U \cap S, S, S', \psi|_S)$ and 
$g_q = g(q(U), V/S, V'/S', \psi_q)$.
If $(U \cap S, S, S', \psi|_S)$, and 
$(q(U), V/S, V'/S', q'\circ \psi \circ q^{-1})$ are good, then
$g_p$ and $g_q$ are respectively
dense subsets of the $B$-orbits 
$B \cdot p(U) \subset Gr(S)$ and  $B \cdot q(U) \subset Gr(V/S)$.
By $B$-equivariance of $p$ and $q$, 
$p^{-1}(g_p) \cap B \cdot U$ and $q^{-1}(g_q) \cap B \cdot U$ 
are both dense subsets of $B \cdot U \subset Gr(V)$.  

Take 
$\tilde U \in p^{-1}(g_p) \cap q^{-1}(g_q)$.  Then 
$\psi|_{p(\tilde U)}: \tilde U \cap S \to S'$ and
$\psi_q|_{q(\tilde U)}: q(\tilde U) \to V'/S'$ are both injective.
By elementary linear algebra, 
$\psi|_{\tilde U}: {\tilde U} \to V'$ is therefore
also injective, as required.
\end{proof}

\subsubsection{Factoring through an intermediate module}
In Sections \ref{sec:igame} and \ref{sec:bgame},
the geometric ideas in Propositions \ref{prop:geommove} and 
\ref{prop:geomsplit} will translate into the combinatorics of the root
game.  Our next proposition is not used, because it is not so 
easy to make combinatorial in its full generality.  However, a special
case of this can be nicely incorporated into the root game for 
Schubert intersection numbers; this appears in Section
\ref{sec:computations}.

Given a quadruple $(U,V,V',\psi)$, let $\tilde B$ be a
group such that $B' \subset \tilde B \subset B$, and let $\tilde V$
be a $\tilde B$ module.  Suppose the map $\psi$ factors as 
$\psi = \psi_2 \circ \psi_1$, where $\psi_1 : V \to \tilde V$
is $\tilde B$ equivariant, and $\psi_2: \tilde V \to V'$ is
$B'$ equivariant.

\begin{proposition}
\label{prop:geommerge}
If $\psi_1|_U: U \to \tilde V$ is injective and 
$(\psi_1(U), \tilde V, V', \psi_2)$ is good, then $(U,V,V', \psi)$
good.
\end{proposition}

\begin{proof}
If there exists $a \in \tilde B$ such that $\psi_2|_{a \cdot \psi_1(U)}$
is injective, then $\psi|_{a \cdot U}$ is also injective.
\end{proof}

\subsection{Questions}
The results of this section (Propositions \ref{prop:geommove},
\ref{prop:geomsplit} and \ref{prop:geommerge})  provide a way of
proving that $\branchproblem \neq 0$, by producing a set of varieties
($B$-orbit closures on Grassmannians $Gr(V)$), and $T$-fixed 
points $U$ on these varieties such that $\psi|_U$ is injective. 
A natural question is whether such a $T$-fixed 
point always exists if $\branchproblem=0$.

This question as stated is somewhat vague, and can be phrased more 
precisely in a couple of different ways.  The most obvious interpretation
is does there exist a suitable $T$-fixed point which can be found using 
only the results of this section?  Less restrictively, one might observe
that successive uses of Proposition \ref{prop:geommove} may not find
all the $T$-fixed points on a $B$-orbit.  If one includes all the
$T$-fixed points in the picture, does a suitable $T$-fixed point
always exists?  If so how does one practically find these other 
$T$-fixed points?

The first formulation of the question 
is essentially asking for a converse to Theorems
\ref{v-thm:win} and \ref{b-thm:win}, and unfortunately the answer is 
in general no (see Section \ref{sec:computations} for further
discussion).  The second 
formulation is open and appears to be a difficult problem.  
In \cite{P2} we show 
that the answer is yes for the multiplication problem in the special case 
where the Schubert classes are pulled back from a Grassmannian.


\section{Root games for Schubert intersection numbers}
\label{sec:igame}

\subsection{Overview of the root game}
The root game combinatorially encodes the geometric notions of 
Section \ref{sec:geom}.
We discuss two versions of the game.
In this section  we will handle the special case of the vanishing 
problem of Schubert intersection numbers on $G'/B'$.
We present the most general version (for the branching
problem of a general inclusion $G' \into G$)
in Section \ref{sec:bgame}.

The former is actually a special case of the latter.  We present
the two formulations separately, since several of the rules become
simplified in the root game for Schubert intersection numbers, and
moreover
it is convenient to encode the data slightly differently for these
two problems.  

The basic overview of the game is the same for both problems.  
The playing
field is a set of squares, which correspond to positive roots of
$G$ or $G'$.  Some of the squares contain {\em tokens}, which get moved
from square to square by the player according to certain rules.
The set of all squares is subdivided into {\em regions}, which limit
the movement of the tokens.  The player alternates between subdividing
the regions further ({\em splitting}), and moving around tokens, in an
attempt to reach a {\em winning position}.

In Section \ref{sec:icombgeom}, we shall see the connection with 
the geometry in Section \ref{sec:degenerate}.  
In short, the positions of the tokens will represent
the $T$-weights of potential degenerations of the subspace $Q \subset \frn$.  
The position of the tokens before and after a move will be the weights 
before and after 
a degeneration of the type in Proposition \ref{prop:geommove}, whereas
the splitting of regions corresponds to the type of subdivision in
Proposition \ref{prop:geomsplit}.  Ultimately the purpose of the game is
to search for a degeneration of $Q$ which will allow us to easily
conclude that $\strconst \neq 0$ from Lemma \ref{b-lem:linalg}; 
these will be the winning positions.

We advise the reader who wishes to skip directly to the more general
branching root game to glance first at examples \ref{ex:a-squares} and 
\ref{ex:b-squares}, which illustrate how the squares arranged for
root systems of types A and B, as this will be essential to 
understanding subsequent examples.


\subsection{Rules of the game}

Recall the problem: given $(\pi_1 \ldots, \pi_s)$, the vanishing
problem is to determine whether $\strconst =
\int_{G'/B'} \omega_{\pi_1} \cdots \omega_{\pi_s} = 0$.
We assume that $\sum \ell(\pi_i) = \dim G'/B'$, otherwise this integral
vanishes for dimensional reasons.

\subsubsection{Data of a position}
The {\bf position} in a root game consists of the following data:
\begin{itemize}
\item A partition of the set of positive roots of $G'$, i.e.
$\regions = \{R_1, \ldots, R_r\}$, such that
$\Delta'_+ = \coprod_{i=1}^s R_i$.  Each $R_i$ is called
a {\bf region}.
\item A list of subsets $\tokens_1, \ldots, \tokens_s$ of the positive 
roots of $G'$, which we call the {\bf arrangement of tokens}.
\end{itemize}

We organise these data as follows.  We draw a set of squares:
the squares correspond to the positive roots of $G'$, and are 
arranged in a sensible way (depending on the type of $G'$.)
The squares are denoted $S_\alpha$, $\alpha \in \Delta'_+$.

\begin{example}
\label{ex:a-squares}
Suppose $G' = SL(n)$.  
Let $x_1, \ldots, x_n$ denote an orthonormal basis for $\RR^n$.
The root system $\Delta' = A_{n-1}$ is
$\{\alpha_{ij} = x_j-x_i\ |\ i \neq j\}$.
The positive roots are those
for which $i<j$.
We can view our squares corresponding to the positive roots
as being arranged inside an $n \times n$ array of squares.
Let $AS_{ij}$ denote the square in position $(i,j)$.
The relevant squares are squares $AS_{ij}$ (the square in position
$(i,j)$), where $1 \leq i<j \leq n$.  Thus the positive root
$\alpha_{ij}$, with $i<j$ is assigned to the square $AS_{ij}$.
                                                                                
$$ \Delta(SL(6))_+ = (A_5)_+ =
\text{
\begin{tabular}{ccccc} \hline
\sqr{\alpha_{12}} &
\sqr{\alpha_{13}} &
\sqr{\alpha_{14}} &
\sqr{\alpha_{15}} &
\sqrend{\alpha_{16}} \\\hline
& \sqr{\alpha_{23}} &
\sqr{\alpha_{24}} &
\sqr{\alpha_{25}} &
\sqrend{\alpha_{26}} \\\cline{2-5}
& & \sqr{\alpha_{34}} &
\sqr{\alpha_{35}} &
\sqrend{\alpha_{36}} \\\cline{3-5}
& & & \sqr{\alpha_{45}} &
\sqrend{\alpha_{46}} \\\cline{4-5}
& & & & \sqrend{\alpha_{56}} \\\cline{5-5}
\end{tabular}
}
$$
\end{example}

Each square may contain one or more {\em tokens}.  
We think of the tokens as physical objects which can be
moved from one square to another.
Each token has a label $k \in \onetos$.  Two tokens with the
same label can never be in the same square. 
We'll call a token labeled $k$ a $k$-token, and write $k \in S_\alpha$
if a $k$-token appears in square $S_\alpha$.
The subsets $\tokens_1, \ldots, \tokens_s \subset \Delta_+'$ are always
defined as:
$$\tokens_k := \{ \alpha \in \Delta'_+
\ |\ \text{the square $S_\alpha$ contains a $k$-token}\}.$$

\subsubsection{Initial position}
In the initial position of the game, 
there is only one region: $\regions = \{\Delta'_+\}$.
The arrangement of tokens is the inversion set for
$\pi_1, \ldots, \pi_s$:
$$\tokens_k = \{\alpha \in \Delta'_+\ |\ \pi_k(\alpha) \in \Delta'_-\}.$$

\begin{figure}[htbp]
  \begin{center}
    \epsfig{file=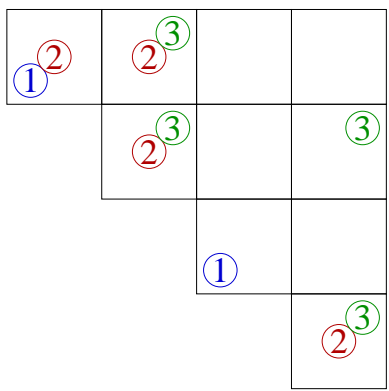,height=2in}
    \caption{Initial position of the game for permutations 21435, 32154, 24153.}    
   \label{v-fig:initialposition}
  \end{center}
\end{figure}

\begin{example}
For $G'=SL(n)$, $\pi_1, \ldots, \pi_s$ are given by permutations of
$1, \ldots, n$.  $\alpha_{ij}$ is an inversion of $\pi(k)$ if and only 
if $\pi(i) > \pi(j)$.  The initial position of the root game for
$G'=SL(5)$, $\pi_1 = 21435$, $\pi_2 = 32154$, $\pi_3 =24153$ 
is given in Figure \ref{v-fig:initialposition}.
\end{example}

From the initial position the player performs a sequence of 
{\em splittings}, which change the set of regions, and {\em moves},
which change the arrangement of tokens.

\subsubsection{Splitting}
\label{sec:igame-split}
Before each move, the player subdivides the regions $R \in \regions$,
according to the following rules.

\begin{definition}
\label{def:idealsubset}
Let $A = \{S_\alpha \ |\  \alpha \in I\}$ be a subset of the squares.
Call $A$ an {\bf ideal subset}\footnote{This is sometimes called an order
ideal for the root poset of $\frn'$.} if $I$ is closed under raising
operations,  i.e. If $\alpha \in I$, then $\alpha' \in I$, whenever
$\alpha'$, and $\alpha'-\alpha$ are both positive roots.
(Equivalently, $A$ is a an ideal subset if and only if
$\{e_\alpha\ |\ \alpha \in I\}$ span an ideal in the Lie algebra $\frn'$.)
\end{definition}

For any ideal subset $A$, we define the operation of 
{\bf splitting along $A$}, as follows: we subdivide each region $R$
into two regions $R \cap A$ and $R \cap A^c$.
(Empty regions produced in this way can be ignored.)  Thus $\regions$
is replaced by 
$$\regions' = \{R_1 \cap A, R_1 \cap A^c,\, R_2 \cap A, R_2 \cap A^c,\, \ldots,
\,R_r \cap A, R_r \cap A^c\}.$$

In principle the player may split along any arbitrary collection of ideal 
subsets between moves; however, this is inadvisable.  The player
{\em should} split along an ideal subset $A$ if and only if the
total number of tokens in the squares of $A$ is
exactly equal to $\#(A)$.  If this condition is followed, each new region 
will always have the property that the number of tokens within the region 
is equal to the number of squares in the region.  
When splitting  is performed with every $A$ satisfying this condition, we 
call the process {\bf splitting maximally}.  No choice is involved in 
splitting maximally.  

\subsubsection{Moving}
After the regions have been split maximally, the player makes a move.  
A move is specified by a triple $[k, \beta, R]$, where $k \in \onetos$
is a choice of token label, $\beta \in \Delta'_+$, and $R \in \regions$ 
is a choice of region.  

To execute the move $[k, \beta, R]$ we change the arrangement of tokens 
as follows.  Find all pairs of squares $S_\alpha, S_\alpha' \in R$ such
that $\alpha'-\alpha = \beta$, and proceeding in order of decreasing
height of $\alpha$, if a $k$-token occurs in the square
$S_\alpha$ but not in $S_{\alpha'}$, move it from the first 
square to the second square.

Using Definition \ref{def:shift} the result of a move can be described 
as follows.  If
$\tokens'_1, \ldots, \tokens'_s$ represents the arrangement
of tokens after the move $[k, \beta, R]$, then for any region
$R' \in \regions$,
$$\tokens'_j \cap R' = 
\begin{cases}
\tokens_j \cap R', &\text{if $R' \neq R$ or $j \neq k$} \\
(\tokens_k \cap R) \shift{R} \beta &\text{otherwise.}
\end{cases}
$$

\subsubsection{Play of the game}
Beginning with the initial position, the player alternates
between splitting maximally (to subdivide the regions), and making a 
move to change the arrangement of tokens.  

\begin{definition}
The game is {\bf won} if at any point there is exactly one token in each
square.
\end{definition}

Observe that a token can only ever move from a square $S_\alpha$ to
$S_{\alpha'}$, where $\alpha'$ is a higher root than $\alpha$.
So, for example, if there are two tokens in the square corresponding
to the highest root, there is no point in proceeding further.  
More generally, 
the game is {\bf lost} if there is an ideal subset $A$ such that the
the total number of tokens in $A$ is more than $\#(A)$.  

\begin{figure}[htbp]
  \begin{center}
    \epsfig{file=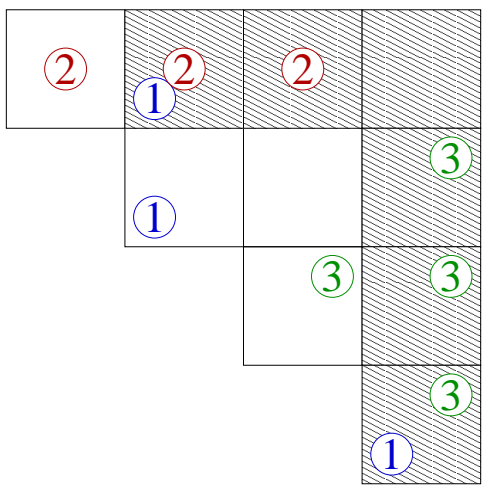,height=2.5in}
    \caption{The initial position for permutations 23154, 41235, 13542
      is a losing position, as there are there are 7 tokens in the 6 
      shaded squares.}
    \label{v-fig:losingexample1}
  \end{center}
\end{figure}

An important special case is when the initial configuration of
tokens is a losing position.  An example of this is shown in Figure
\ref{v-fig:losingexample1}.

\begin{definition}
If the game is lost before the first move is made, we say the game
is {\bf doomed}.
\end{definition}

Note that, while a doomed game cannot be won, it is not the case 
that all games which cannot be won are doomed (as seen in Section 
\ref{sec:convcount}).

\subsubsection{Vanishing and non-vanishing criteria}
From games which are doomed, and games which can be won, we obtain vanishing
and non-vanishing criteria respectively.

\begin{theorem}
\label{v-thm:lose}
If the root game corresponding to $(\pi_1, \ldots, \pi_s)$
 is doomed, then $\strconst = 0$.
\end{theorem}

\begin{theorem}
\label{v-thm:win}
If the root game corresponding to $(\pi_1, \ldots, \pi_s)$ can be
won, then $\strconst \geq 1$.
\end{theorem}

\begin{remark}
It is also possible
to play the game, omitting the splitting stage.  This simplifies
the combinatorics considerably, and Theorem
\ref{v-thm:win} still holds.  However, as mentioned already, it is
always advisable to split maximally between moves.  It is easy to show 
that if the game can be won by omitting the splitting step, it can be
still be won while including it (c.f. Section \ref{sec:splitting}).
\end{remark}

\begin{remark}
In the case where the game is doomed as a result of an ideal subset $A$
which is maximal (i.e. $A$ consists of all squares except for a single
$S_\alpha$, where $\alpha$ is a simple root), Theorem \ref{v-thm:lose}
reduces to the DC-triviality vanishing condition
in \cite{Kn}.
\end{remark}

Theorems \ref{v-thm:lose} and \ref{v-thm:win} cast a large net over
the set of all Schubert problems, and capture a huge number of them. 
It is not hard to see, for instance, that the 
probability of finding a non-doomed game at random for $SL(n)$ tends to $0$ 
as $n \to \infty$.  Still there is a small gap: in
general, not being able to win the game does not provide any 
information.  However, in a number special cases, we have been able to
show that the converse of Theorem \ref{v-thm:win} holds.  These are
discussed in Section \ref{sec:vremarks}.
                                                                                

\subsection{Relating the combinatorics and geometry}
\label{sec:icombgeom}
Given an arrangement of tokens, $\tokens_1, \ldots, \tokens_s$, and
a region $R \in \regions$, we associate the following linear spaces.
\begin{itemize}
\item A $B'$-module $V'$, a subquotient of $\frn'$, such that $R$ is 
the set of distinct $T'$-weights of $V'$.
\item For each $i=1, \ldots, s$, a linear subspace $U_i \subset V'$ 
such that $\tokens_i \cap R$ is the set of $T'$ weights of $U_i$.
\end{itemize}
Put $V = V' \oplus \cdots \oplus V'$ ($s$-summands), and 
$U = U_1 \oplus \cdots \oplus U_s \subset V$.  We have a $B'$-equivariant
map $\phi: V \to V'$ given by $\phi(v_1, \ldots, v_s) = v_1 + \cdots + v_s$.
Thus we have a quadruple $(U,V,V',\phi)$, as in Section \ref{sec:degenerate}.

At any position in the game, there is one such quadruple for every region.
Recall the notion of a {\em good} quadruple from Section \ref{sec:degenerate}.
In our current context, the quadruple $(U,V,V',\phi)$ is good if and only if
there exist $a_1, \ldots, a_s \in N$ such that $a_1 \cdot U_1 \oplus \cdots
\oplus a_s \cdot U_s = V'$.  A position is a winning position if and only
if $U_1 \oplus \cdots \oplus U_s = V'$ for every region.  In particular,
the quadruples associate to winning positions are good.

To prove Theorem
\ref{v-thm:win}, we show that if the quadruple associated to every region
is good, then the quadruple associated to initial position is good.
This is true essentially because the moves and splittings combinatorially 
encode 
the geometric ideas in Propositions \ref{prop:geommove} and 
\ref{prop:geomsplit}.
A move in the root game changes $U$ (for a single region)
to a new subspace $U^1 \subset V'$ in exactly the manner prescribed
by Proposition \ref{prop:geommove}. Thus, if the new
quadruple $(U^1, V, V', \phi)$ is good, then old quadruple $(U,V,V',\phi)$
must be good too.  Similarly a splitting changes the set of quadruples
following Proposition \ref{prop:geomsplit}.  Hence from a winning position 
(where all regions correspond to good quadruples),
we can backtrack all the way to the start the game and deduce that the
initial position is good.

Now, the statement that the initial position is good is precisely 
condition 2 of 
Corollary \ref{cor:v-linalg}.  Thus, by Corollary \ref{cor:v-linalg}, 
the initial position is good if and only if $\strconst \neq 0$.

Many of the details have been omitted here.  A more precise account 
of relationship between the geometry and combinatorics is given in the 
proofs in Section \ref{sec:proofs}.


\subsection{Examples}

\subsubsection{Games which can be won}
\begin{figure}[htbp]
  \begin{center}
    \epsfig{file=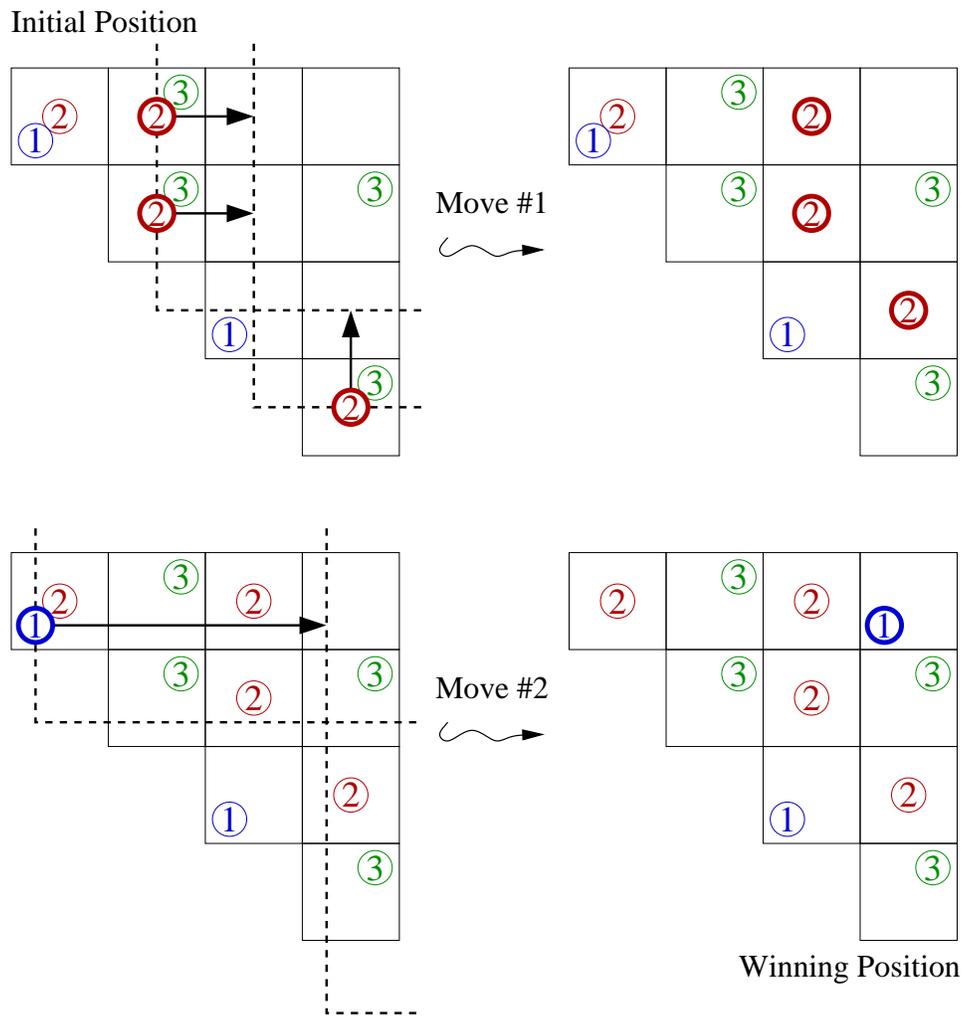,width=5in}
    \caption{Moves $[2, \alpha_{34}]$ and $[1, \alpha_{25}]$ are 
      applied to the initial
      position in Figure \ref{v-fig:initialposition}.}
    \label{v-fig:winningexample}
  \end{center}
\end{figure}

\begin{figure}[htbp]
  \begin{center}
    \epsfig{file=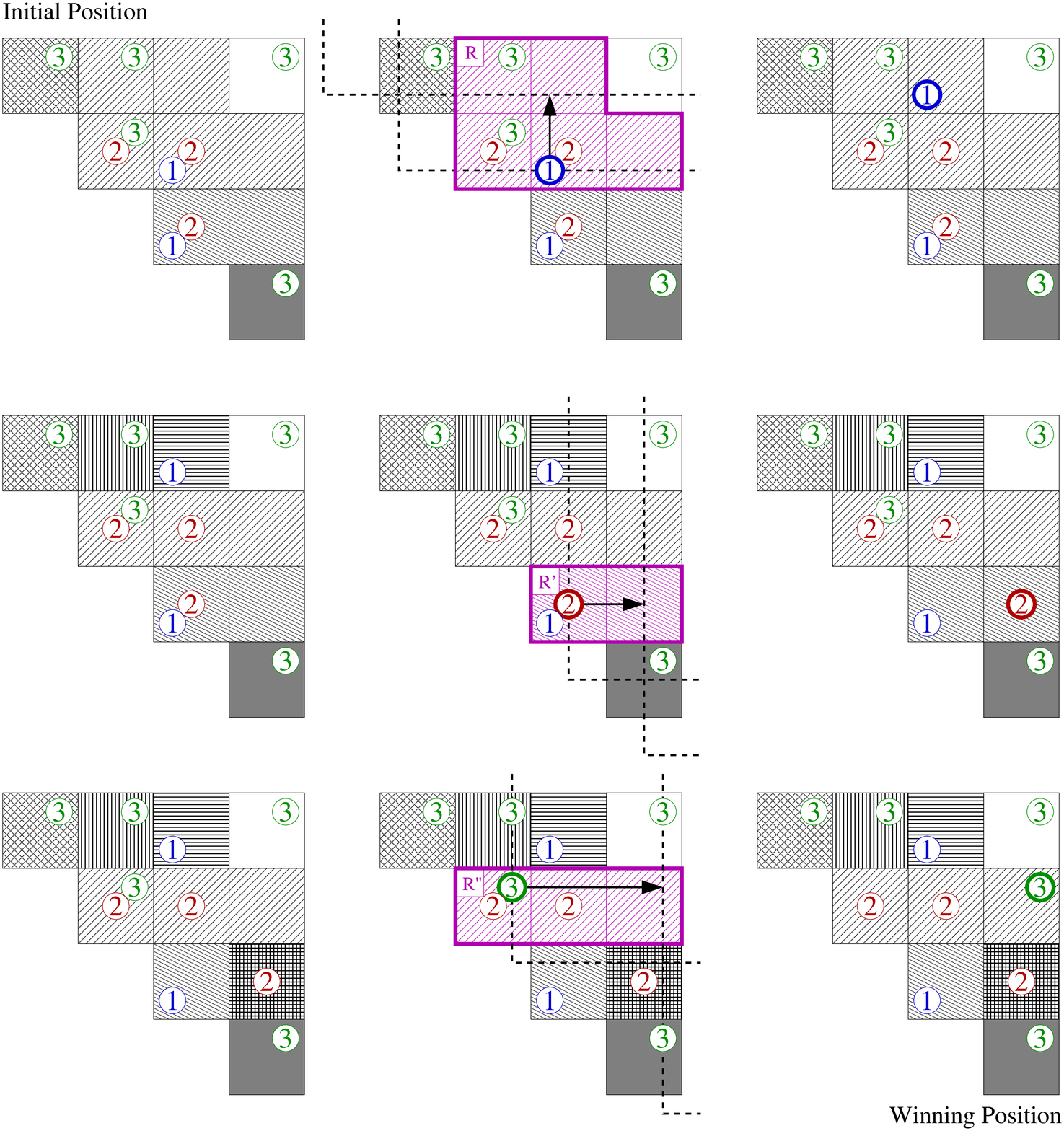,width=5.8in}
    \caption{The general game, played out for permutations 13425, 41325,
14352.  The moves, shown in the centre column, are:  $[1,\alpha_{12},R]$,
$[2,\alpha_{45},R']$, and finally $[3,\alpha_{35},R'']$.
The left column shows the state before the move, in which the set of
squares is maximally divided into regions. The right column shows the
state immediately after the move, before further subdividing.}
    \label{v-fig:generalexample}
  \end{center}
\end{figure}

In type $A$, for any
fixed $\beta=\alpha_{ij}$, the possible squares involved in a move 
corresponding to $\beta$ are
$ \{S_{\alpha_{kl}}\ |\ \text{$k=i$ or $l=i$}\}
\cup \{S_{\alpha_{kl}}\ |\ \text{$k=j$ or $l=j$}\}
$.
These squares lie on two reflected lines
which meet at the square $S_\beta$ (shown as dotted lines in
Figures \ref{v-fig:winningexample} and \ref{v-fig:generalexample}).  
The tokens move strictly horizontally or vertically from one reflected 
line the other.

Figure \ref{v-fig:winningexample} shows a sequence of moves in a
game which has been played without the splitting step, to better illustrate
the movement of the tokens.  The initial position
is $(21435, 32154, 24153)$, from Figure \ref{v-fig:initialposition}.
The sequence of moves leads to a winning position.

Figure \ref{v-fig:generalexample} gives an example of a sequence
of moves with maximal splitting in between moves.  Again the sequence
of moves leads to a winning position.

\subsubsection{Converses and counterexamples}
\label{sec:convcount}
\begin{figure}[htbp]
  \begin{center}
    \epsfig{file=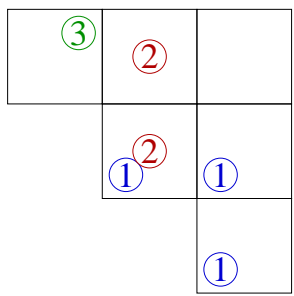,height=1.5in}
    \caption{The permutations $\pi_1=1432$, $\pi_2=2314$, $\pi_3=2134$ are
a counterexample to the converse of Theorem \ref{v-thm:lose}.  The game is
not doomed, though $\triplestrconst=0$.  All other $SL(4)$ counterexamples are
similar to this one.}
    \label{v-fig:losingconverse}
  \end{center}
\end{figure}

The converse of Theorem \ref{v-thm:lose} is certainly not true.  The first
counterexamples in $SL(n)$ occur for $n=4$.  See Figure
\ref{v-fig:losingconverse}.

\begin{figure}[htbp]
  \begin{center}
    \epsfig{file=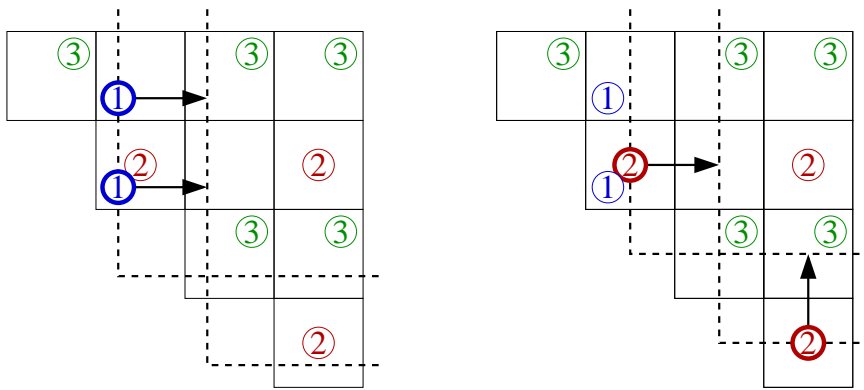,height=2in}
    \caption{The permutations 23145, 14253, 41523 give a counterexample to the
converse of Theorem \ref{v-thm:win}, if the game is played without splitting.
If we split maximally first, either of the moves shown (restricted to the
appropriate region) will win the game.}
    \label{v-fig:wimpyconverse}
  \end{center}
\end{figure}

If the root game is played without splitting, the converse of 
Theorem \ref{v-thm:win} is not true.
The first counterexamples in $SL(n)$ occur for $n=5$.
Figure \ref{v-fig:wimpyconverse} shows 
the initial position of the game for the permutations 23145, 14253, 41523.
There is only one square with 2 tokens, and one empty square.  Without 
splitting, any effort to rectify this imbalance winds up moving 
more than just one token.
However, with splitting, the game can be won.

\subsubsection{Other types}
                                                                                
We now describe a `sensible' way to arrange the squares in types
$B$ and $D$ ($G'= SO(n)$).  A similar arrangement to the type $B$
arrangement can be used for type $C$ ($G' = Sp(2n)$).

In both examples $x_1, \ldots, x_n$ is an orthonormal basis
for $\RR^n$.

\begin{example} 
\label{ex:d-squares}
If $G' =SO(2n)$, the root system $\Delta'=D_n$ is 
$$\{ (-1)^{\varepsilon} x_i + (-1)^{\delta} x_j\ |\ i \neq j \}.$$
The positive roots are of two types:
$$\{\beta{ij} = x_j -x_i\ |\ i < j\} \cup
\{\beta'_{ij} = x_j+x_i\ |\ i < j\}$$
We arrange these the squares inside a $2n \times n$ array 
(denoted $DS_{ij}$)
as follows: the root $\beta_{ij}$ corresponds to the square 
$DS_{n+i,j}$;  the root $\beta'_{ij}$ corresponds to the square
$DS_{n+1-i,j}$.
$$ \Delta(SO(10))_+ =(D_5)_+ =
\text{
\begin{tabular}{cccc}  \cline{4-4}
&&& \sqrend{\beta'_{45}} \\\cline{3-4}
& & \sqr{\beta'_{34}} &
\sqrend{\beta'_{35}} \\\cline{2-4}
& \sqr{\beta'_{23}} &
\sqr{\beta'_{24}} &
\sqrend{\beta'_{25}} \\\hline
\sqr{\beta'_{12}} &
\sqr{\beta'_{13}} &
\sqr{\beta'_{14}} &
\sqrend{\beta'_{15}} \\\hline
\sqr{\beta_{12}} &
\sqr{\beta_{13}} &
\sqr{\beta_{14}} &
\sqrend{\beta_{15}} \\\hline
& \sqr{\beta_{23}} &
\sqr{\beta_{24}} &
\sqrend{\beta_{25}} \\\cline{2-4}
& & \sqr{\beta_{34}} &
\sqrend{\beta_{35}} \\\cline{3-4}
&&& \sqrend{\beta_{45}} \\\cline{4-4}
\end{tabular}
}$$

\end{example}

\begin{example} 
\label{ex:b-squares}
If $G' =SO(2n+1)$, the root system $\Delta'=B_n$ is 
$$\{ (-1)^{\varepsilon} x_i + (-1)^{\delta} x_j\ |\ i \neq j \}
\cup \{ \pm x_i\}.$$
The positive roots are of three types:
$$\Delta'_+=
\{\gamma_{ij} = x_j -x_i\ |\ i < j\} \cup
\{\gamma'_{ij} = x_j+x_i\ |\ i < j\} \cup
\{\gamma^\circ_j = x_j\}.$$
We arrange the squares inside a $(2n+1) \times n$ array of squares
(denoted $BS_{ij}$)
as follows: the root $\gamma_{ij}$ corresponds to the square 
$BS_{n+1+i,j}$;  the root $\gamma'_{ij}$ corresponds to the square
$BS_{n+1-i,j}$;  the root $\gamma^\circ_j$ corresponds to the square
$BS_{n+1,j}$.

$$\Delta(SO(9))_+ = (B_4)_+ =
\text{
\begin{tabular}{cccc}  \cline{4-4}
&&& \sqrend{\gamma'_{34}} \\\cline{3-4}
& & \sqr{\gamma'_{23}} &
\sqrend{\gamma'_{24}} \\\cline{2-4}
& \sqr{\gamma'_{12}} &
\sqr{\gamma'_{13}} &
\sqrend{\gamma'_{14}} \\\hline
\sqr{\gamma^\circ_1} &
\sqr{\gamma^\circ_2} &
\sqr{\gamma^\circ_3} &
\sqrend{\gamma^\circ_4} \\\hline
& \sqr{\gamma_{12}} &
\sqr{\gamma_{13}} &
\sqrend{\gamma_{14}} \\\cline{2-4}
& & \sqr{\gamma_{23}} &
\sqrend{\gamma_{24}} \\\cline{3-4}
&&& \sqrend{\gamma_{34}} \\\cline{4-4}
\end{tabular}
}$$
\end{example}

\begin{figure}[htbp]
  \begin{center}
    \epsfig{file=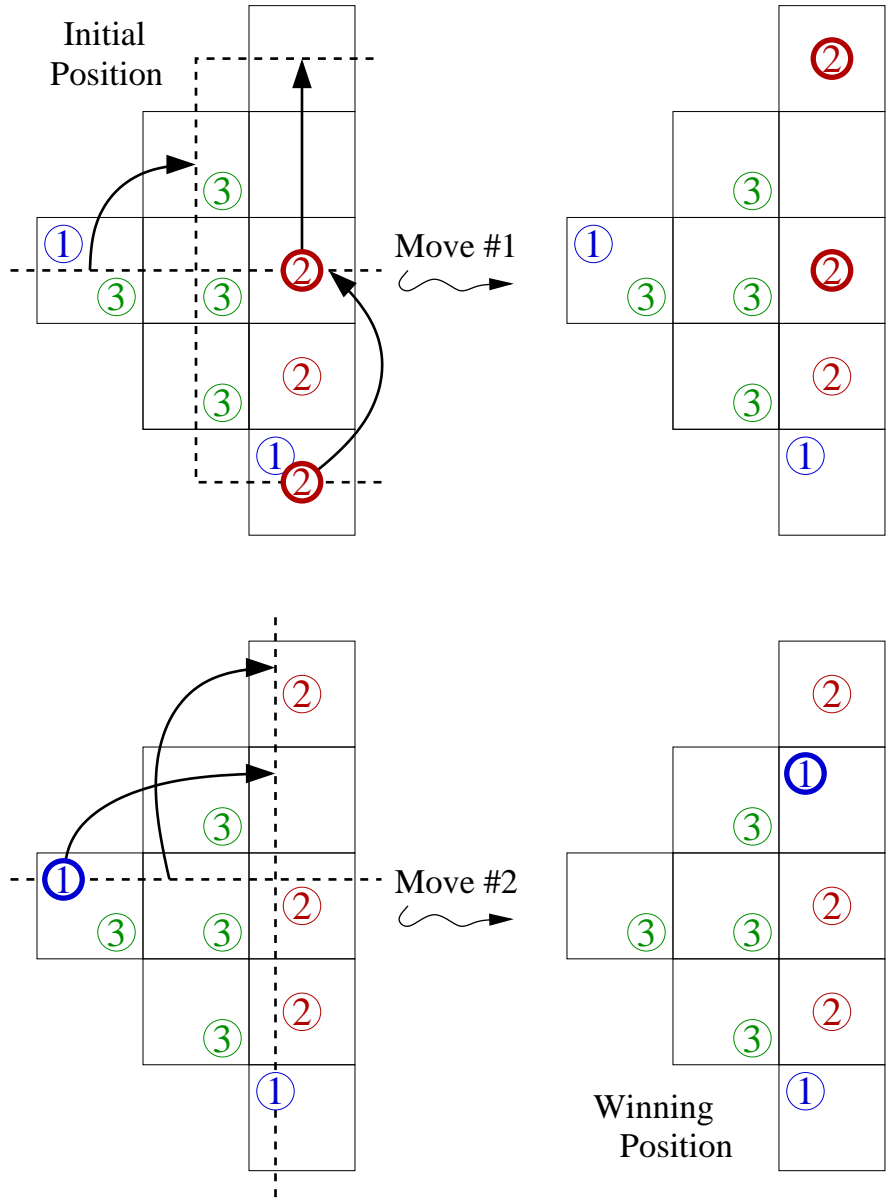,width=4.6in}
    \caption{A simple game for $SO(7,\CC)$.  In this example,
$\pi_1 = \bar{1}32$, $\pi_2 = 231$, $\pi_3=\bar{1}\bar{2}3$.
The root $\beta$ which is used in each move is the crossing point
of the dotted lines.}
    \label{v-fig:b3example}
  \end{center}
\end{figure}

Figure \ref{v-fig:b3example} gives an example of a root 
game for $SO(7,\CC)$.
Here, an element of Weyl group $W=C_2^3 \rtimes S_3$ can be
represented by a permutation $a_1 a_2 a_3$ of $1 2 3$, where each
symbol is either decorated with a bar or not.  This permutation
acts on $\RR^3$ by the matrix whose $i^{th}$ row is $x_{a_i}$ if
$i$ is unbarred, and $-x_{a_i}$ if $i$ is barred.
                                                                                
Arrows in Figure \ref{v-fig:b3example}
are included not only for all tokens that move, but for all pairs of
roots $\alpha, \alpha'$, whose difference is $\beta$.  
Since the game
can be won using the moves shown, for
$\pi_1 = \bar{1}32$, $\pi_2 = 231$, $\pi_3=\bar{1}\bar{2}3$
we have $\triplestrconst \geq 1$.

\subsection{Remarks}
\label{sec:vremarks}

\subsubsection{Splitting}
\label{sec:splitting}
In the rules of the root game, we are told exactly when to
split regions: we split 
along an ideal subset $A$ if and only if the number of tokens in 
$A$ equals $\#(A)$.  
However, it turns out that this condition is never used in the proof. 
Thus, in theory, the rules could be relaxed so that 
the player has 
the option to split regions along {\em any} ideal subset $A$ between moves.
That said, we will now sketch a proof that it is never advantageous to 
the player to exercise this freedom.

Suppose the rules tell us not to split along $A$.  If we do split
along $A$ there will be too many tokens in one region.  Since regions
can never be rejoined once they are split, the game cannot be won.
On the other hand, 
suppose the rules tell us to split along $A$, and the player chooses
not to.  Of any move that is made subsequently, one of the following two
things must be true: {\em either} the same arrangement of tokens could have 
been reached (possibly using multiple moves) if we had split along $A$, 
{\em or}
the move caused the game to be lost.

The ability to determine, a priori when splitting is advantageous,
relies on the fact that we have assumed 
$\sum_{i=1}^s \ell(\pi_i) = \dim G'/B'$.

\subsubsection{Products which are not top degree}
\label{sec:nontopdegree}
The root game can be adapted to analyse
the non-vanishing of a product
of Schubert classes, whether or not the product is not of top degree.  
This greater generality is handled by the root game for branching 
Schubert calculus, so we will not discuss it at length.
However only two minor modifications to the rules are
required.  First, we must change the winning condition
to read ``the game is won if there is at most one token in each square'', 
rather than 
``... exactly one token in each square''.  (The losing and 
doomed conditions remain as stated previously.)  Second, we must
remove the rule forcing us to split maximally between moves, and instead
have splitting be at the player's discretion, as 
discussed in Section \ref{sec:splitting}.

\subsubsection{Relationship with the Bruhat order}

For products of only two Schubert
classes, the converse of Theorem \ref{v-thm:win} holds:
being able to win the root game both necessary and
sufficient for non-vanishing.
In this case,
the non-vanishing of the product is determined precisely by the
Bruhat order.  That is,
$\omega_{\pi_1} \cdot \omega_{\pi_2} \neq 0$ if and only
if $\pi_1 \leq w_0 \pi_2$ in the Bruhat order.  

In the case where the product is top degree, i.e. 
$\ell(\pi_1) + \ell(\pi_2) = \dim G/B$, the fact that we can win the
root game is a triviality: we have
$\pi_1 \leq w_0 \pi_2$ if and only if $\pi_1 = w_0 \pi_2$, in which
case the set of squares containing a $1$-token is the complement
of the set of squares containing a $2$-token.  Thus the initial
position of the game is already a winning position.

Less trivial is the case when $\pi_1 < w_0 \pi_2$.  Since the product of
the classes is not top degree, we must use the revised notion of
winning position (see Section \ref{sec:nontopdegree}).  Nevertheless,
we have the following theorem.

\begin{theorem}
\label{v-thm:bruhat}
$\omega_{\pi_1} \cdot \omega_{\pi_2} \neq 0$
if and only if it is possible to win the root game
corresponding to $(\pi_1, \pi_2)$.
\end{theorem}

A detailed proof of this result is given in the author's doctoral thesis
\cite{P2}.

\subsubsection{Converses and computations}
\label{sec:computations}
It would be quite surprising and remarkable if the converse of Theorem
\ref{v-thm:win} were true in any generality.  So far, for $SL(n)$, 
the converse has deftly eluded any counterexamples.  In fact the 
converse of Theorem \ref{v-thm:win} has been affirmed by an 
exhaustive computer search for $SL(n)$ for $n \leq 7$.  The converse
of Theorem \ref{v-thm:win} (with $s=3$)
has also been verified for the exceptional group $G_2$, as well
as for $SO(5)$ and $SO(7)$.  (The next
smallest exceptional group, $F_4$, is
unfortunately beyond our computational abilities at the moment.)

One special case where the converse of Theorem \ref{v-thm:win} is
true is when the classes $\omega_{\pi_i}$ are pulled back from 
a Grassmannian in an appropriate way.  We prove this result
in \cite{P1}.  Another special case is
Theorem \ref{v-thm:bruhat}, which tells us that the converse is
true for products of only two Schubert classes.

For the groups $SO(n)$, $n \geq 8$, the converse of Theorem 
\ref{v-thm:win} is in fact false.  For $SO(8)$, we represent
an element of the Weyl group
$W=C_2^3 \rtimes S_4$ 
by a permutation $a_1 a_2 a_3 a_4$ of $0 1 2 3$, where 
each symbol, except $0$, is either decorated with a bar or not ($0$ is 
always unbarred).  This permutation acts on $\RR^4$ by the matrix 
whose $i^{th}$ row is $x_{a_i}$ if $i$ is unbarred, and $-x_{a_i}$ if 
$i$ is barred.  The two counterexamples to the converse of Theorem 
\ref{v-thm:win} for $SO(8)$ are listed below.
\begin{center}
\begin{tabular}{|c|c|c|}
\hline
$\pi_1$ & $\pi_2$ & $\pi_3$ \\
\hline \hline
$0\bar{1}32$ & $0\bar{2}31$ &  $03\bar{2}1$ \\\hline
$03\bar{1}2$ & $0\bar{2}31$ &  $0\bar{2}31$ \\\hline
\end{tabular}
\end{center}
The problem that arises in these examples
is that although there exists a $T$-fixed point $(U_1,U_2,U_3)$ on
$Gr(V')^3$ with 
$U_1 \oplus U_2 \oplus U_3 = V'$, which is a degeneration of 
$(a_1 \cdot P_1, a_2 \cdot P_2, a_3 \cdot P_3)$,
the moves of the game fail to find it.
We are not aware of any examples in which $\strconst \geq 1$ but where there 
are no suitable $T$-fixed points on any of the relevant varieties.
It therefore seems it would be desirable to be able to describe a larger 
set of moves---moves which, starting from a $T$-fixed point on $Gr(V')^s$ 
can reach all the other $T$-fixed points in its $(N')^s$-orbit closure.

With $s=3$, a restriction one might wish to make to the root game
is to allow only moves involving tokens labeled $1$ and $2$.  This 
restriction seems appropriate when viewing Schubert calculus as taking 
products in
cohomology rather than intersection numbers.  
Under this weakening, Theorem \ref{v-thm:win}
remains true (obviously), but the converse is already false for $SL(n)$.  
There are no examples of this for $n \leq 5$; however, for $n=6$ there
are a total of four such examples:
\begin{center}
\begin{tabular}{|c|c|c|}
\hline
$\pi_1$ & $\pi_2$ & $\pi_3$ \\
\hline \hline
145326 & 321564 &  315264 \\\hline
154326 & 312564 &  315264 \\\hline
514326 & 152364 &  135264 \\\hline
154236 & 312654 &  315264 \\\hline
\end{tabular}
\end{center}
It is possible
dispose of these $SL(6)$ counterexamples, by introducing new moves 
geometrically based on Proposition \ref{prop:geommerge}. 
One such move (called a {\em merge}) 
is the following: select a region $R$, and a pair of 
token labels: $k_1 \neq k_2$, with the property that there is no 
square in $R$ which contains both a $k_1$-token and a $k_2$-token.  
Then replace every $k_2$-token with a $k_1$-token in the same square.
If we introduce merges into the root game, the aforementioned 
counterexamples disappear.  Moreover, it follows (though we omit the
proof here) that the converse of 
Theorem \ref{v-thm:win} (with merges included)
is true for $SL(6)$ for any $s$.


\section{Root games for branching Schubert calculus}
\label{sec:bgame}

We now describe the root game for the more general branching problem.  
We will see that although the setup is slightly different, this game 
specialises to the root game for vanishing of Schubert intersection numbers.

The main differences we shall see are the following:
\begin{enumerate}
\item Squares correspond to positive roots of $G$ rather than $G'$.
\item Tokens do not have labels.
\item The winning condition is defined in terms of the map 
$\hat \phi: \Delta_+ \to \Delta'_+ \cup \zeroset$
(which we calculate explicitly in a number of examples).
\item Splitting is more complicated---it also involves $\hat \phi$.
\end{enumerate}


\subsection{The map $\hat \phi$}
Recall the definition of the map $\hat \phi$ from Section 
\ref{sec:conventions}.
The rules of the game heavily involve $\hat \phi$; thus
before proceeding further, we compute this map in a number
of important examples.

\begin{example} 
If $G' \into G_1, \ldots, G' \into G_s$, then 
$G' \into G = G_1 \times \cdots \times G_s$
via the diagonal map.  
Let $\hat \phi_i: \Delta(G_i)_+ \to \Delta' \cup \zeroset$ denote the
map on root systems for $G' \into G_i$.
The positive roots of $G$ are 
$\Delta_+ = \Delta(G_1)_+ \sqcup \dots \sqcup \Delta(G_s)_+$, and 
$$\hat \phi: \Delta_+ \to \Delta'_+ \cup \zeroset$$
is simply given by
$ \hat \phi(\alpha) = \phi_i(\alpha)$ if $\alpha \in \Delta(G_i)_+$.
\end{example}
In particular, if $G'=G_1 = \cdots = G_s$, then each $\hat \phi_i$ is
just the identity map.  Thus this example allows us to deal 
with the vanishing problem for multiplication of Schubert calculus.

\begin{example} 
If $G' = SL(k) \into G = SL(n)$ is the inclusion
$$
A \mapsto
\left (
\begin{matrix}
A & 0 \\
0 & I_{n-k}
\end{matrix}
\right )
$$
then 
$$\hat \phi(\alpha_{ij}) =
\begin{cases}
 \alpha_{ij} \in \Delta', &\text{if $j \leq k$} \\
 0, &\text{otherwise}
\end{cases}
$$
In this example $i^*$ is easily described as an operation on Schubert
polynomials \cite{LS}.
If $\omega_\pi$ is represented by a Schubert polynomial
in variables $x_1, \ldots, x_n$ then $\branchproblem$ is given by
setting $x_{k+1} = \cdots = x_n = 0$.
\end{example}

We now consider the inclusion of $SO(n, \CC) \hookrightarrow G = SL(n)$.
We begin with the case where $n$ is even.  
Let $R$ denote the $n/2 \times n/2$ matrix with $1$ on the antidiagonal, 
and $0$ everywhere else.  We take as a (compact) maximal torus of 
$SO(n, \RR)$ the subgroup
$$ T'_\RR = \big \{
\left(
\begin{matrix}
A & BR \\
-RB & RAR 
\end{matrix}
\right ) \in GL(n, \RR)
\ \big |\ \text{$A$, $B$ are diagonal} \}.$$
The complexification $T'$ of $T'_\RR$ is a complex maximal torus in 
$SO(n,\CC)$.

Since
the standard maximal torus of $SO(n)$ is not a subgroup of the standard 
(diagonal) maximal torus of $SL(n)$, we use a suitable conjugate subgroup
of $SO(n)$.  Let
$$U = \left (
\begin{matrix}
I & iR \\
iR & I
\end{matrix}
\right ).$$

\begin{figure}
\begin{center}
\begin{tabular}{ccccccc} \hline
\sqr{\beta_{34}} &
\sqr{\beta_{24}} &
\sqr{\beta_{14}} &
\sqr{\beta'_{14}} &
\sqr{\beta'_{24}} &
\sqr{\beta'_{34}} &
\sqrend{} \\\hline 
& 
\sqr{\beta_{23}} &
\sqr{\beta_{13}} &
\sqr{\beta'_{13}} &
\sqr{\beta'_{23}} &
\sqr{} &
\sqrend{\beta'_{34}} \\\cline{2-7}
&& 
\sqr{\beta_{12}} &
\sqr{\beta'_{12}} &
\sqr{} &
\sqr{\beta'_{23}} &
\sqrend{\beta'_{24}} \\\cline{3-7}
& & &
\sqr{} &
\sqr{\beta'_{12}} &
\sqr{\beta'_{13}} &
\sqrend{\beta'_{14}}  \\\cline{4-7}
& & & &
\sqr{\beta_{12}} &
\sqr{\beta_{13}} &
\sqrend{\beta_{14}}  \\\cline{5-7}
& & & & &
\sqr{\beta_{23}} &
\sqrend{\beta_{24}} \\\cline{6-7}
& & & & &&
\sqrend{\beta_{34}} \\\cline{7-7}
\end{tabular}
\end{center}
\caption{The map $\hat \phi: \Delta(SL(8))_+ \to \Delta(SO(8))_+ \cup \zeroset$.
The root $\hat \phi(\alpha)$ is written in the square corresponding to
$\alpha$.  Empty squares are mapped to $0$.}
\label{b-fig:da-rootmap}
\end{figure}

\begin{example} 
\label{ex:da-rootmap}
Let $G' = U\,SO(n)\,U^{-1} \into G = SL(n)$, where $n=2m$.  One
can easily verify that a maximal torus of $G'$, is the set of 
invertible diagonal matrices
$$U\,T'\,U^{-1} = \big \{\lambda =
\left(
\begin{matrix} 
\lambda_{m} & & & & & \\
& \ddots & & & & \\
& & \lambda_1 & & & \\
& & & \lambda_1^{-1} & & \\
& & & & \ddots &  \\
& & & & & \lambda_{m}^{-1}
\end{matrix}
\right)
\in SL(n) \big \}.
$$
The Lie algebra $\frg' = \{(a_{ij}) \in \frsl(n)\ |\ a_{ij}=-a_{n+1-i\ j}\}$,
is the set of $n \times n$ matrices which are skew symmetric about the
antidiagonal.  And $\frn'$ is simply the set of upper triangular matrices
in $\frg'$.  

Let $E_{ij}$ denote the matrix with a $1$ in the $i,j$ position, and 
$0$ everywhere else. 
We see that for $i<j$,
$$\lambda E_{ij} \lambda^{-1} =
\begin{cases}
\lambda_{m+1-i}\lambda_{j-m} E_{ij}, 
&\text{if $i+j>n+1$, $i \leq m$} \\
\lambda_{m+1-i}\lambda_{j-m} E_{ij}, 
&\text{if $i+j \leq n$, $j>m$} \\
\lambda_{m+1-i}\lambda_{m+1-j}^{-1} E_{ij}, 
&\text{if $j \leq m$} \\
\lambda_{i-m}^{-1}\lambda_{j-m} E_{ij}, 
&\text{if $i > m$} \\
0, &\text{if $i+j=n+1 $} \\
\end{cases}
$$

Thus $\hat \phi$ is given by
$$\hat \phi(\alpha_{ij}) = 
\begin{cases}
\beta'_{{m+1-i},{j-m}},
&\text{if $i+j>n+1$, $i \leq m$} \\
\beta'_{{j-m},{m+1-i}},
&\text{if $i+j \leq n$, $j>m$} \\
\beta_{{m+1-j},{m+1-i}},
&\text{if $j \leq m$} \\
\beta_{{i-m},{j-m}},
&\text{if $i > m$} \\
0, &\text{if $i+j=n+1 $} \\
\end{cases}
$$

In terms of the arrangement of squares (described in Examples
\ref{ex:a-squares} and \ref{ex:d-squares}), the map $\hat \phi$ 
is symmetrical about the antidiagonal, with the antidiagonal itself 
mapping to $0$.  Moreover, below the antidiagonal (i.e. for $i+j>n+1$), 
we simply have $\hat \phi(AS_{ij}) = DS_{i-m,j}$.  See Figure
\ref{b-fig:da-rootmap}.
\end{example}

The analysis for $n$ odd is very similar.

\begin{example} 
\label{ex:ba-rootmap}
Let $G' \into G = SL(n)$,
where $n=2m-1$ and $G'$ is an appropriately chosen conjugate of $SO(n)$.  

As in the case where $n$ is even, 
$\frg' = \{(a_{ij}) \in \frsl(n)\ |\ a_{ij}=-a_{n+1-i,j}$,
is the set of $n \times n$ matrices which are skew symmetric about the
antidiagonal, and $\frn' = \frb \cap \frg'$.  The map $\hat \phi$ is given by
$$\hat \phi(\alpha_{ij}) = 
\begin{cases}
\gamma^\circ_i,
&\text{if $j=m$} \\
\gamma^\circ_j,
&\text{if $i=m$} \\
\gamma'_{{m-i},{j-m}},
&\text{if $i+j>n+1$, $i<m$} \\
\gamma'_{{j-m},{m-i}},
&\text{if $i+j \leq n$, $j>m$} \\
\gamma_{{m-j},{m-i}},
&\text{if $j < m$} \\
\gamma_{{i-m},{j-m}},
&\text{if $i > m$} \\
0, &\text{if $i+j=n+1 $} \\
\end{cases}
$$

More simply, in terms of the arrangement of squares (see Examples
\ref{ex:a-squares} and \ref{ex:b-squares}), we have that
$\hat \phi$ is symmetrical about the antidiagonal, and identically
zero on the antidiagonal.  
Below the antidiagonal $\hat \phi(AS_{ij}) = BS_{i,j-m}$.  See
Figure \ref{b-fig:ba-rootmap}.
\end{example}

\begin{figure}
\begin{center}
\begin{tabular}{cccccc} \hline
\sqr{\gamma_{23}} &
\sqr{\gamma_{13}} &
\sqr{\gamma^\circ_3} &
\sqr{\gamma'_{13}} &
\sqr{\gamma'_{23}} &
\sqrend{} \\\hline
& \sqr{\gamma_{12}} &
\sqr{\gamma^\circ_2} &
\sqr{\gamma'_{12}} &
\sqr{} &
\sqrend{\gamma'_{23}} \\\cline{2-6}
&& \sqr{\gamma^\circ_2} &
\sqr{} &
\sqr{\gamma'_{12}} &
\sqrend{\gamma'_{13}} \\\cline{3-6}
& & &
\sqr{\gamma^\circ_1} &
\sqr{\gamma^\circ_2} &
\sqrend{\gamma^\circ_3} \\\cline{4-6}
& & & & \sqr{\gamma_{12}} &
\sqrend{\gamma_{13}} \\\cline{5-6}
& & & &  &\sqrend{\gamma_{23}} \\\cline{6-6}
\end{tabular}
\end{center}
\caption{The map $\hat \phi: \Delta(SL(7))_+ \to \Delta(SO(7))_+ \cup \zeroset$.
The root $\hat \phi(\alpha)$ is written in the square corresponding to
$\alpha$.  Empty squares are mapped to $0$.}
\label{b-fig:ba-rootmap}
\end{figure}

\begin{figure}[htbp]
  \begin{center}
    \epsfig{file=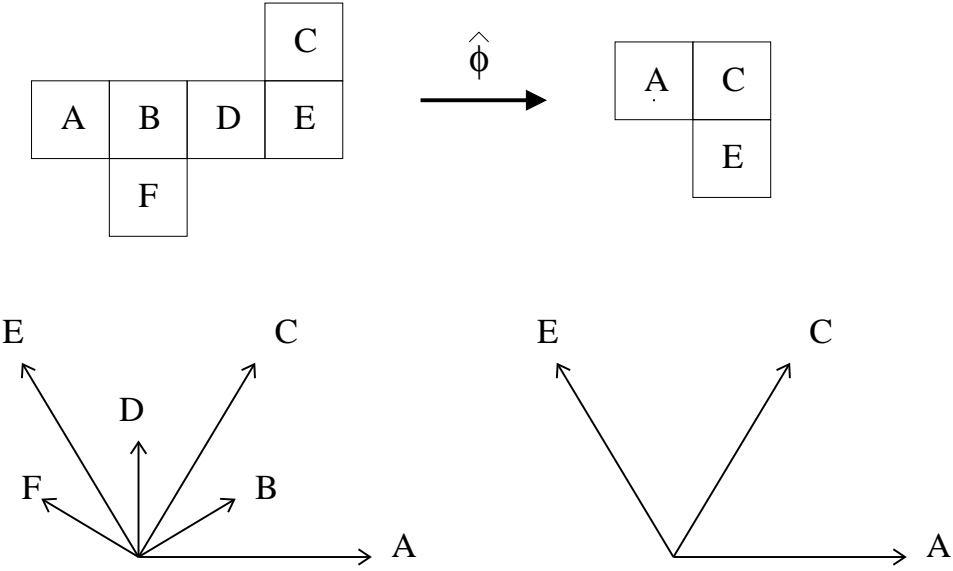,width=5in}
    \caption{The map $\hat \phi: (G_2)_+ \to (A_2)_+ \cup \zeroset$, and the 
corresponding arrangement of squares.}
    \label{b-fig:ag-rootmap}
  \end{center}
\end{figure}

\begin{example}
\label{ex:ag-rootmap}
Let $G$ be the complex form of $G_2$, and $G' = SL(3)$.  The map
$i: G' \into G$ is defined on the level of roots: $A_2$ includes into
$G_2$ as the long roots.  Since $SL(3)$ is simply connected, this
defines a homomorphism on the Lie groups 
(and this map is an inclusion).
The map 
$\hat \phi: (G_2)_+ \to (A_2)_+ \cup \zeroset$ is therefore 
$$
\hat \phi(\alpha) =
\begin{cases}
0, \qquad &\text{if $\alpha$ is a short root of $G_2$} \\
\alpha, \qquad &\text{if $\alpha$ is a long root of $G_2$.}
\end{cases}
$$
We arrange the squares of $G$ in a linear fashion, with the short simple
root at the bottom, and the long simple root on the left.
The map $\hat \phi$ and the arrangement of squares for $G_2$ are
both illustrated in Figure \ref{b-fig:ag-rootmap}.
\end{example}

%


\subsection{Rules of the game}
\subsubsection{Data of a position}

The {\bf position} in a root game consists of the following data:
\begin{itemize}
\item A partition of the set of positive roots of $G$, i.e.
$\regions = \{R_1, \ldots, R_r\}$, such that
$\Delta_+ = \coprod_{i=1}^s R_i$.  Each $R_i$ is called
a {\bf region}.
\item A subset $\tokens$ of the positive roots of $G$, which we call the
{\bf arrangement of tokens}.
\end{itemize}

We visualise this information by drawing a {\em square} $S_\alpha$
for each positive root $\alpha \in \Delta_+$, 
and placing a {\em token} in $S_\alpha$ if $\alpha \in T$.
As before, we arrange the squares in a sensible manner depending
on the type of $G$ (see Examples \ref{ex:a-squares}, 
\ref{ex:d-squares} and \ref{ex:b-squares}).

The regions are just sets of the squares.  As such, if $R$ is
a region, we will sometimes write $S_\alpha \in R$ rather than 
$\alpha \in R$.

\subsubsection{Initial configuration}
\label{sec:movesplit}

The game always begins with a single region ($ \regions = \{\Delta_+\}$), 
which contains all the squares.  The initial arrangement of tokens is the
inversion set of $\pi$, i.e.
$$\tokens = \{ \alpha \in \Delta_+\ |\ \pi \cdot \alpha \in \Delta_-\}.$$

\begin{example}
\label{ex:initialpos}
If $G = SL(5) \times SO(5)$, $\pi = (23154, r_{\gamma^\circ_1})$
where $r_{\gamma^\circ_1})$ is the reflection in the simple root
$\gamma^\circ_1$.  Then the initial position is as shown below:

\begin{center}
 \epsfig{file=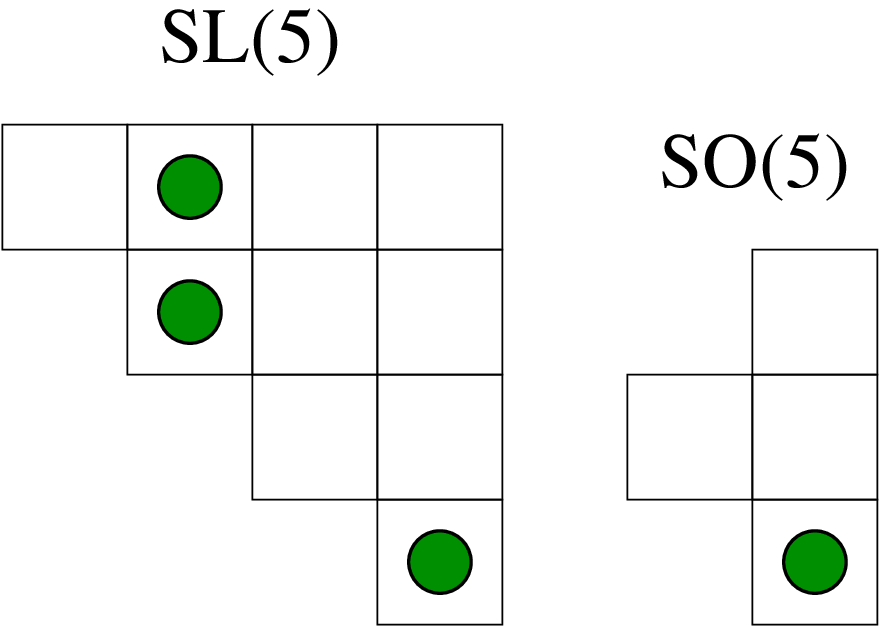,height=1in}
\end{center}

\end{example}

\subsubsection{Splitting}
We define {\em splitting along $A$} as in Section \ref{sec:igame-split}:
if $A \subset \Delta_+$, splitting 
$\regions = \{R_1, \ldots, R_r\}$ along $A$ produces
$$\regions' = \{R_1 \cap A, R_1 \cap A^c,\, R_2 \cap A, R_2 \cap A^c,\, \ldots,
\,R_r \cap A, R_r \cap A^c\}.$$
(Empty regions have no effect on the game, thus we may discard any copies 
of the empty set produced in this way.)

The subsets
$A \subset \Delta$ which can be legally used for splitting are called
\emph{splitting subsets}.  

\begin{definition}
Let $A \subset \Delta_+$ be a subset of the positive roots of $G$.  We
call $A$ a {\bf splitting subset} if $A$ 
is an ideal subset (c.f. Definition \ref{def:idealsubset}), and 
$\hat \phi(A^c) \cap \hat \phi(A) \subset \zeroset$.
\end{definition}

\begin{example}
For $SO(n) \into SL(n)$, a set $A \subset \Delta_+$ is an ideal
subset if for every square $S$ in $A$, $A$ contains all
squares above and to the right of $S$.
$A$ is a splitting subset if it is an ideal subset which is
symmetrical about the antidiagonal.
\end{example}

\subsubsection{Moves}
A {\bf move} is specified by a pair $[\beta, R]$, where
$\beta \in \Delta_+$, and $R \in \regions$ is a choice of region.  
To execute
the move, we find all pairs of squares $S_\alpha, S_\alpha' \in R$ such
that $\alpha'-\alpha = \beta$.  We then order the relevant $S_\alpha$
according to the height of the root $\alpha$.
Proceeding in order of decreasing
height of $\alpha$, if a token appears
in the square $S_\alpha$ but not in $S_{\alpha'}$, move the token up
from the first square to the second square.

Equivalently, using Definition \ref{def:shift} the result of a move can 
be described as follows.  If
$\tokens'$ represents the arrangement
of tokens after the move $[\beta, R]$, then for any region
$R' \in \regions$,
$$\tokens' \cap R' = 
\begin{cases}
\tokens \cap R', &\text{if $R' \neq R$} \\
(\tokens \cap R) \shift{R} \beta. &\text{otherwise}
\end{cases}
$$

\subsubsection{Play of the game}
Beginning with the initial configuration, the player performs
a sequence of moves and splittings.  Moves and splittings may
be performed in any order.  Splitting along $A$ is permissible
whenever $A$ is a splitting subset.

\begin{definition}
The game is {\bf won} if the arrangement of tokens $\tokens$ is injective
(c.f. Definition \ref{b-def:injective}).
\end{definition}

\begin{remark}
Splitting along $A$ is permissible whenever
$A$ is a splitting subset.  However, when $\ell(\pi) = \dim G'/B'$, one can
determine a priori whether splitting will help us win the root game.  
It turns out that 
if $\ell(\pi) = \dim G'/B'$, the splitting is advantageous if and only if
$\#(\tokens \cap A) = \#(\hat \phi(A) \setminus \zeroset)$.  The argument
is analogous to the one given in Section \ref{sec:splitting}.
\end{remark}

From certain positions it may be impossible for victory to be attained.
In particular, if $A$ is an ideal subset, then any
token which begins its move in $A$ must remain in $A$.  Thus
$\#(\tokens \cap A)$ can never decrease over a sequence of moves.  
Suppose then, at some position in the game, there is an ideal subset $A$
such that $\#(\tokens \cap A) > \#(\hat \phi(A) \setminus \zeroset)$.  Then
$\tokens$ is not injective, and will never be injective; thus
the game cannot be won.  In such a position, we declare the game to be 
{\bf lost}.

The situation when the game is lost before any moves are made, is
particularly important.

\begin{definition}
The game is {\bf doomed} if it is lost in the initial token 
arrangement.
\end{definition}

\subsubsection{Vanishing and non-vanishing criteria}
Vanishing and non-vanishing criteria arise from games which are
doomed, and games which can be won.  Games which are lost, or
simply cannot be won provide no information.

\begin{theorem}
\label{b-thm:lose}
If the game is doomed, then $\branchproblem = 0$.
\end{theorem}

\begin{theorem}
\label{b-thm:win}
If the game can be won, then $\branchproblem \neq 0$.
\end{theorem}

Theorems \ref{b-thm:lose} and \ref{b-thm:win} 
specialise to Theorems \ref{v-thm:lose} and \ref{v-thm:win},
taking $i:G' \into G = G' \times \cdots \times G'$ to be the diagonal
inclusion (see Section \ref{sec:specialise}).


\subsection{Examples}

\subsubsection{A corollary of Theorem \ref{b-thm:lose}}

\begin{example}
Let $G = SL(n)$ and $G' = SO(n)$.  Let 
$\pi:\{1, \ldots, n\} \to \{1, \ldots, n\} \in S_n$.  
If $\pi(n) < \pi(1)$ then $\branchproblem =0$.
\end{example}

\begin{proof}
To see this, observe that $A = \{\alpha_{1n}\}$ is an ideal subset,
whose image under $\hat \phi$ is $\zeroset$.  Thus 
$\#(\hat \phi(A) \setminus \zeroset) = 0$.  If $\pi_n < \pi_1$,
then $\alpha_{1n} \in \tokens$, so $\#(\tokens \cap A) = 1$ and
the game is doomed.
\end{proof}

\subsubsection{Games which can be won}

\begin{example}
If $\frg'$ a is $T$-invariant subspace of $\frg$, and 
$\hat \phi^{-1}(\zeroset)$ is an ideal subset, then the initial
position is a winning position if and only if the game is not
doomed, giving a simple necessary and sufficient condition for
$\branchproblem =0$.
Unfortunately this only occurs when the Dynkin diagram of $G'$ is
obtained by deleting some of the vertices of $G$'s Dynkin 
diagram.  Some common examples include $SL(k) \into SL(n)$, 
$SO(2k+1) \into SO(2n+1)$ and $SO(2k) \into SO(2n)$, for $k < n$.
\end{example}

\begin{example}
Let $G = SL(5) \times SO(5)$, $\pi = (23154, r_{\gamma^\circ_1})$.
The initial position is shown in Example \ref{ex:initialpos}.  
We can win the game with one move, and no splittings.  The move 
corresponds
to the root $\gamma^\circ_2 \in (B_2)_+$.  This causes the
token on the $SO(5)$ part to move from $\gamma^\circ_1$ to
$\gamma'_{12}$.
\begin{center}
 \epsfig{file=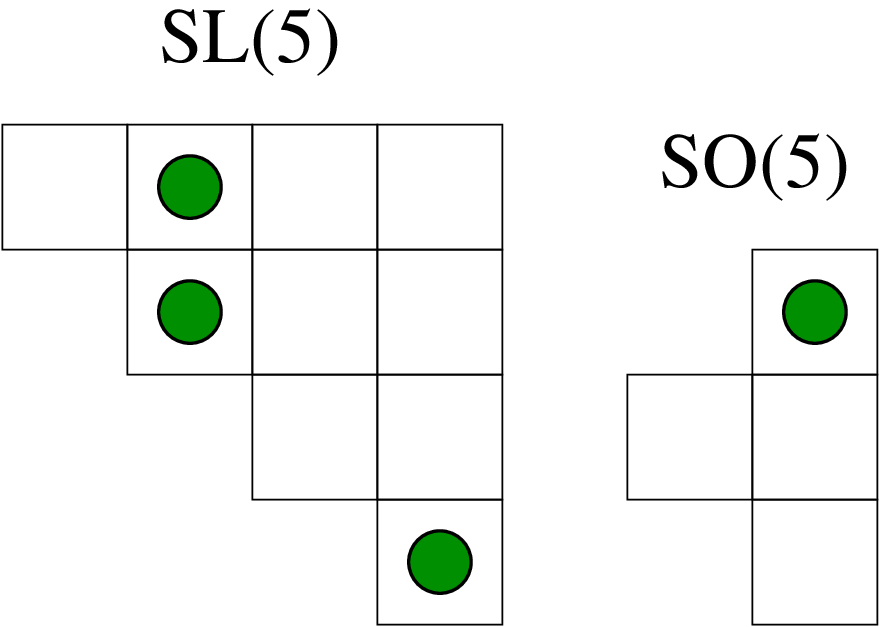,height=1in}
\end{center}
To see that this is a winning position, we fold the $SL(5)$ picture 
along the antidiagonal (this is $\hat \phi: (A_4)_+ \to (B_2)_+$).
\begin{center}
 \epsfig{file=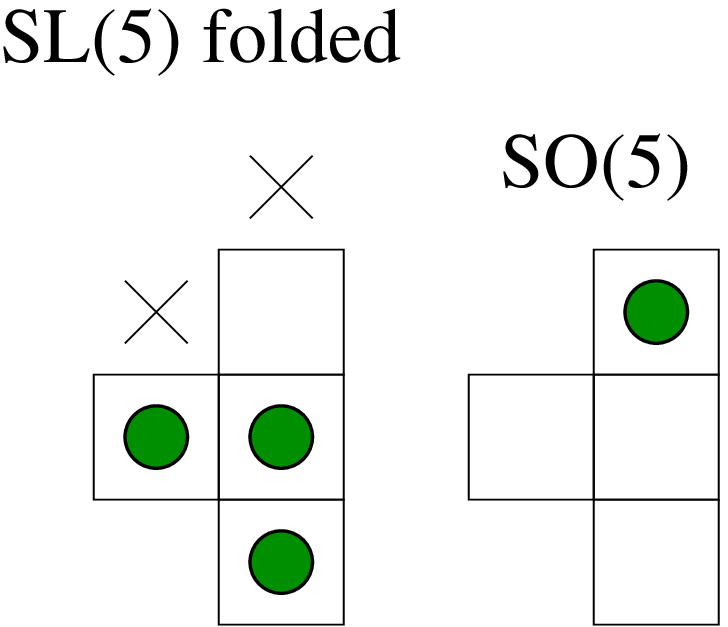,height=1in}
\end{center}
(The {\rm `$\times$'}s denote the diagonal of the folding map.)
We then superimpose the two $(B_2)_+$ pictures which this folding produces
(this is $\hat \phi : (B_2)_+ \times (B_2)_+ \to (B_2)_+$).
Since no tokens overlap in this process, or
appear on the diagonal of the folding map (= $\hat \phi^{-1}(\zeroset)$),
this is a winning position.
\end{example}

\begin{figure}[htbp]
  \begin{center}
    \epsfig{file=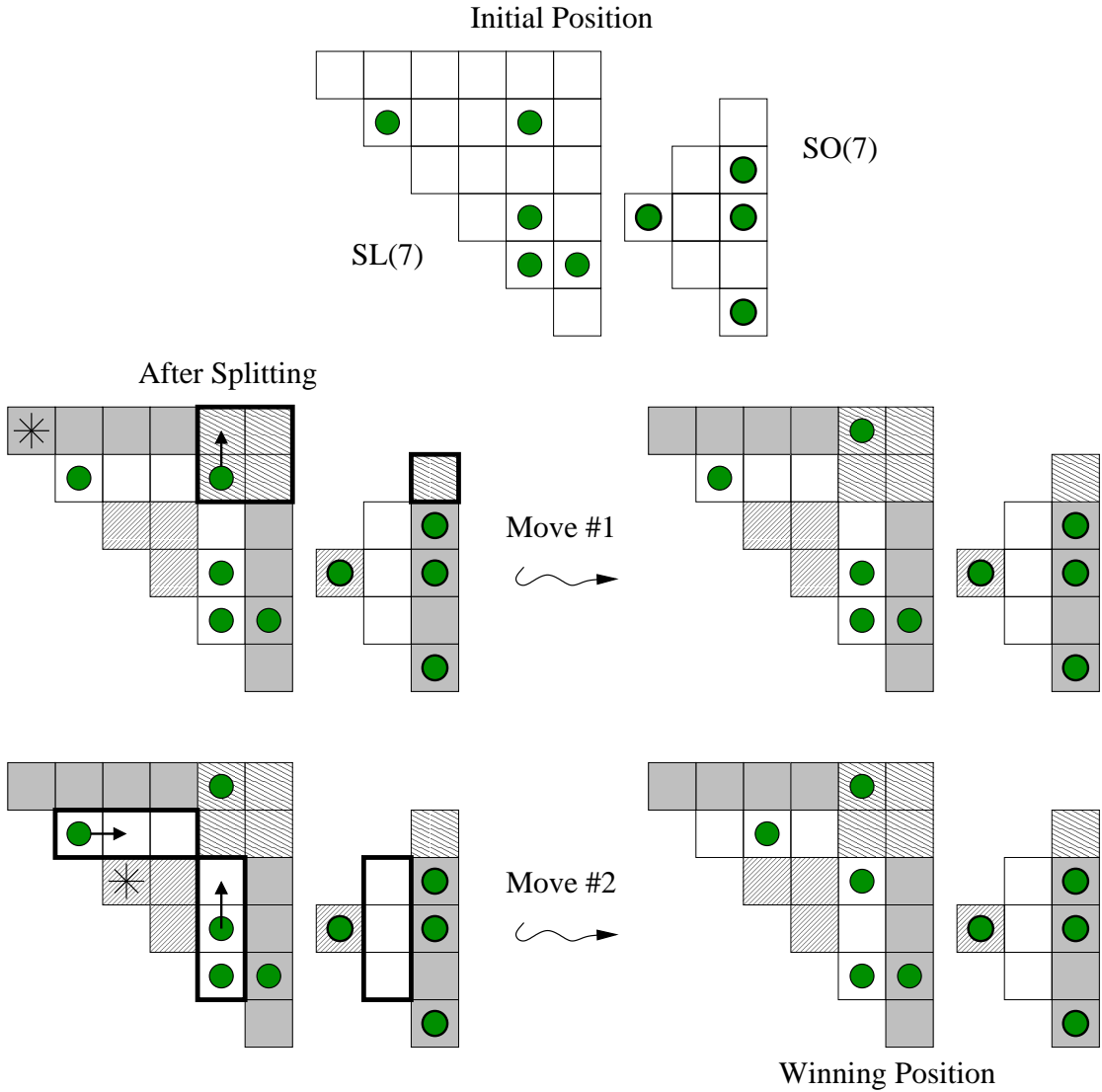,width=5.8in}
    \caption{A sequence of moves in the root game for
$SO(7) \into SL(7) \times SO(7)$, $\pi = (1425736, \bar{2}3\bar{1})$.
The bold outline indicates which region is being used in each move, and
the $*$ indicates which root is being used.}
    \label{b-fig:winexample1}
  \end{center}
\end{figure}

\begin{example}
Let $G = SL(7) \times SO(7)$, $G' = SO(7)$, $\pi = (1425736, \bar{2}3\bar{1})$,
where $\bar{2}3\bar{1}$ is the $SO(7)$ Weyl group element represented
by the matrix 
$$\left (
\begin{matrix}
0 & -1 & 0 \\
0 & 0 & 1 \\
-1 & 0 & 0 
\end{matrix}
\right).
$$
Figure \ref{b-fig:winexample1} shows a sequence of splittings and
moves lead to a winning position.  Squares belonging to the same
region are similarly shaded.  For each move, the relevant region
is outlined, and the relevant root is indicated by an asterisk in
the corresponding square.
\end{example}

\begin{figure}[htbp]
  \begin{center}
    \epsfig{file=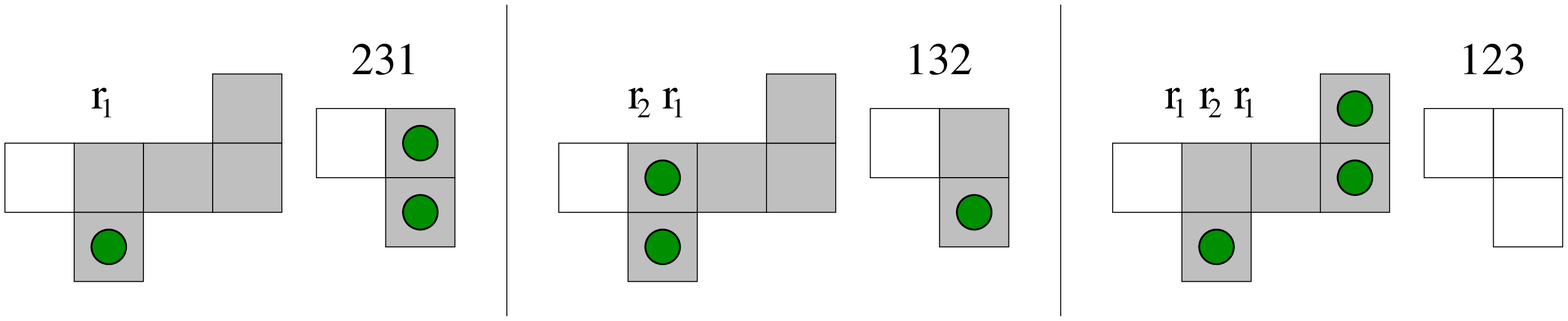,width=5in}
    \caption{The $3$ games which are doomed for 
$SL(3) \into G_2 \times SL(3)$. The shaded squares indicate a minimal
ideal subset $A$ for which 
$\#(\tokens \cap A) > \#(\hat \phi(A) \setminus \zeroset)$.}
    \label{b-fig:ag-doomed}
  \end{center}
\end{figure}
\begin{figure}[htbp]
  \begin{center}
    \epsfig{file=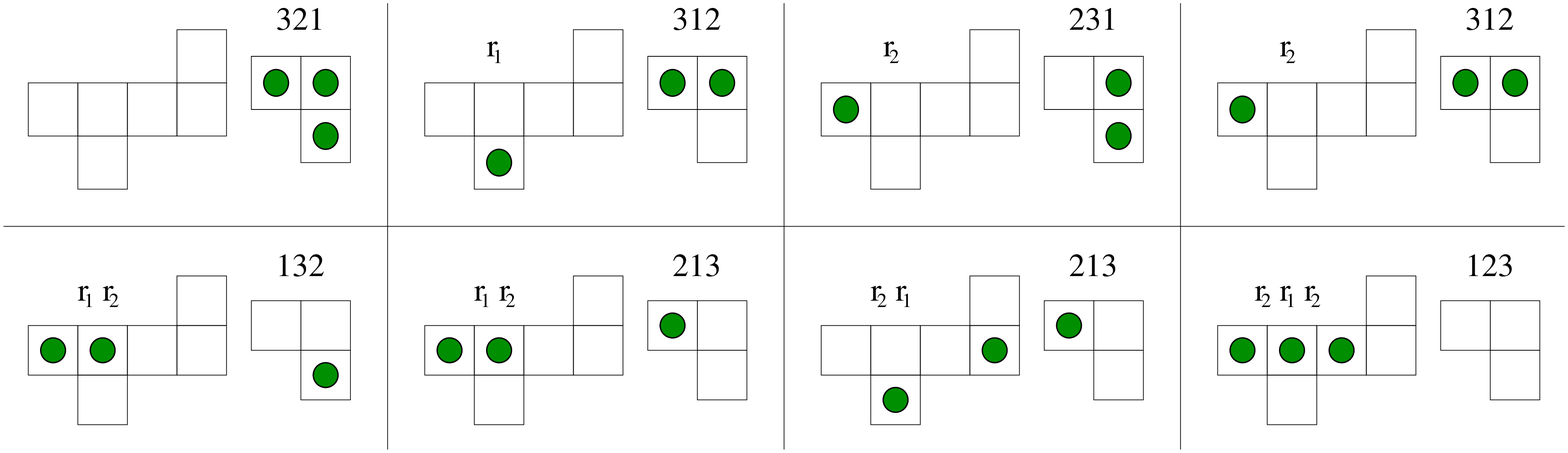,width=5.8in}
    \caption{The $8$ games which are not doomed for 
$SL(3) \into G_2 \times SL(3)$. Each of these games can be won.}
    \label{b-fig:ag-nondoomed}
  \end{center}
\end{figure}
\begin{figure}[htbp]
  \begin{center}
    \epsfig{file=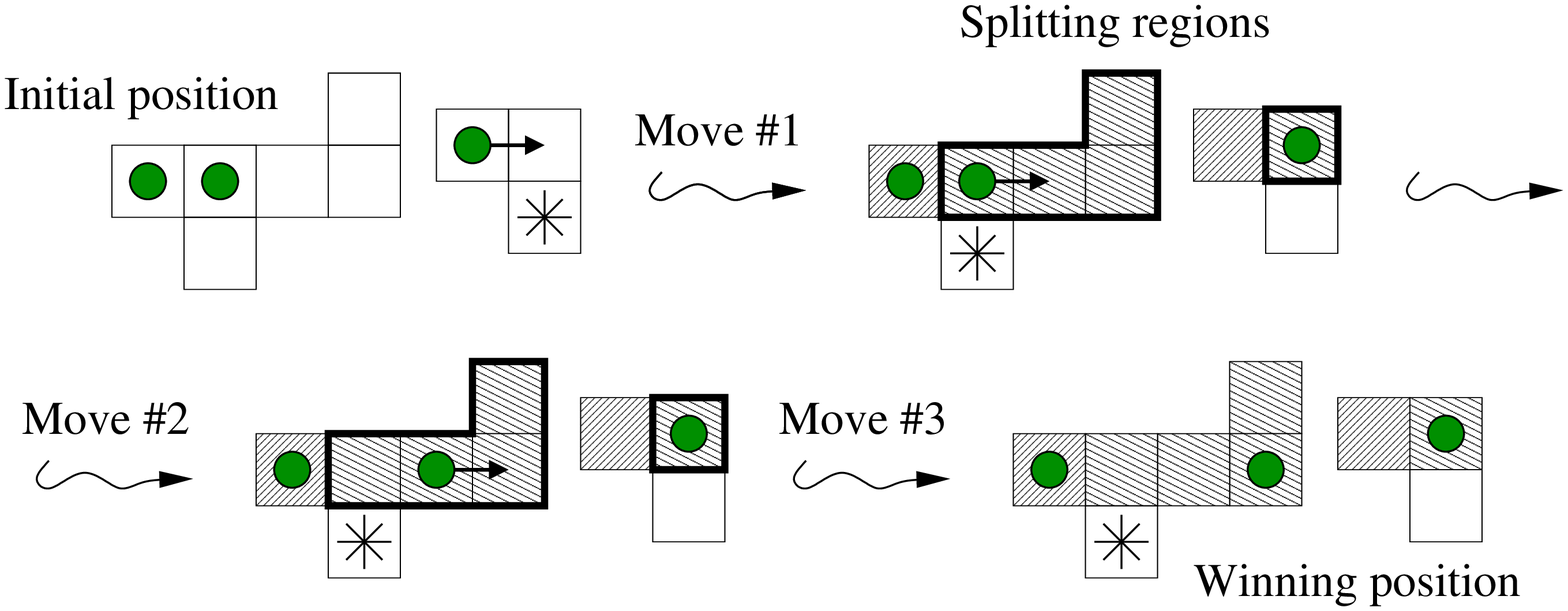,width=5.5in}
    \caption{A sequence of moves in the root game for
$SL(3) \into G_2 \times SL(3)$, $\pi = (r_1 r_2, 213)$.
After the first move, we split into three regions, indicated by the
different shading of squares.
The bold outline indicates which region is being used in each move, and
the $*$ indicates which root is being used.}
    \label{b-fig:winexample2}
  \end{center}
\end{figure}

\begin{example}
Let $G = G_2 \times SL(3)$, and $G' =SL(3)$ including diagonally,
where $SL(3) \into G_2$ is as described in Example \ref{ex:ag-rootmap}.
We consider all possible $\pi \in W$, with $\ell(\pi) = 3 = \dim G'/B'$.
There are such $11$ such $\pi$ in total.  Of these, $3$ associated
games are doomed.
These are $\pi = (r_1, 231)$, $\pi = (r_2 r_1, 132)$, and 
$\pi = (r_1 r_2 r_1, 123)$, where $r_1$ and $r_2$ represent reflections
in the short and long simple roots respectively.  
These are shown in Figure \ref{b-fig:ag-doomed}.
The remaining $8$ games are shown in Figure \ref{b-fig:ag-nondoomed}.
One can check that each of these can be won.  Figure \ref{b-fig:winexample2}
shows a sequence of moves from the initial position of one of these games,
$\pi = (r_1 r_2, 213)$, to a winning position.  Thus the root game
gives a complete answer to the vanishing problem for branching
$SL(3) \into G_2$.
\end{example}

%


\subsection{Specialisation to Schubert intersection numbers}
\label{sec:specialise}

In order to avoid proving Theorems \ref{v-thm:lose} and \ref{v-thm:win}
directly, we show that the formulation of the root game
for vanishing of Schubert intersection numbers
is in fact just a special
case of the more general root game for branching.  In the
interest of brevity, we'll call these two formulations of the
root game the {\em I-game} and the {\em B-game} respectively.

As has been already discussed, the correspondence comes from from putting
$i: G' \into G = G' \times \cdots \times G'$, the diagonal inclusion.

\subsubsection{Token labels versus squares}
In the I-game, squares correspond to positive roots of $G'$, whereas tokens
are labelled $1, \ldots, s$.  In the B-game, the squares correspond
to positive roots of $G$ and the tokens are unlabelled.  The equivalence
of the two is seen
from the fact that the $\Delta_+$ is a disjoint
union of $s$ copies of $\Delta'_+$.  The token label
in the I-game indicates which copy of $\Delta'_+$is being used
for the corresponding token in the B-game.

\subsubsection{Winning condition and splitting}
The map $\hat \phi: \Delta_+ \to \Delta'_+ \cup \zeroset$, is given
by superimposing all the copies of $\Delta'_+$.  Since $\Delta'_+$
corresponds to the set of squares in the I-game, the injectivity
of $\hat \phi|_\tokens$ corresponds to having at most one
token in each square.  Since we assumed that
$\sum \ell(\pi_i) = \dim G'/B'$, this is the same as exactly one
token in each square.

From this description of $\hat \phi$ it is also easy to see that
splitting subsets of $\Delta_+$ are in one to one correspondence
with ideal subsets of $\Delta'_+$.


\subsection{Proofs}
\label{sec:proofs}

\subsubsection{Proof of the vanishing criterion}
The vanishing criterion (Theorem \ref{b-thm:lose}) is a combinatorial
reinterpretation of Lemma \ref{lem:geomlose}.

\begin{proof}[Proof of Theorem \ref{b-thm:lose}]
At the outset, $\{e_\alpha\ |\ \alpha \in \tokens\}$ is a basis for the
space $Q$.  If the game is doomed then there is an ideal subset
$A$ such that $\#(\tokens \cap A) > \#(\hat \phi(A) \setminus \zeroset).$
Let $S \subset \frn$ be the ideal generated by 
$\{e_\alpha\ |\ \alpha \in A\}$.  We have
$$\dim (Q \cap S) = \#(\tokens \cap A) > 
\#(\hat \phi(A) \setminus \zeroset)= \dim \phi(S).$$
By Lemma \ref{lem:geomlose}, we conclude that $\branchproblem = 0$.
\end{proof}

\subsubsection{Proof of the non-vanishing criterion}
The non-vanishing criterion (Theorem \ref{b-thm:win}) is essentially
combinatorially encoding the geometric ideas in Propositions
\ref{prop:geommove} and \ref{prop:geomsplit}.

\begin{proof}[Proof of Theorem \ref{b-thm:win}]
If $R \subset \Delta_+$, we let 
$\langle e_\alpha\ |\ \alpha \in R \rangle_B$
denote the $B$-submodule of $\frn$ generated by all $e_\alpha$,
$\alpha \in R$.

Let $(\regions= \{R_1, \ldots, R_r\}, \tokens)$ be a position in the 
root game.  We associate to this position the following geometric
data:

\begin{itemize}
\item $B$-modules $V_i$, $i=1, \ldots, r$.
Put $V_i^0 = \langle e_\alpha\ |\ \alpha \in R_i \rangle_B$, and
$V_i^1 
= \langle e_\alpha\ |\ \alpha \in \Gamma(V_i^0) \setminus R_i \rangle_B$.
Then we define $V_i$ to be the quotient $V_i^0/V_i^1$.  Note that
$V_i$ is a $B$-module, and a subquotient of $\frn$, with weights
$\Gamma(V_i) = R_i$.

\item $B'$-modules $V'_i$, $i=1, \ldots, r$, defined as 
$V'_i = \phi(V_i^0)/\phi(V_i^1)$.

\item $B'$-equivariant maps $\phi_i:V_i \to V'_i$, induced
from $\phi|_{V_i}$.

\item Subspaces $U_i \subset V_i$.  $U_i$ is defined to be the 
$T$-invariant subspace of $V_i$ with weights 
$\Gamma(U_i) = \tokens \cap R_i$.
\end{itemize}

Thus, for each region $R_i$, $i \in \{1, \ldots, r\}$ we have a
quadruple $(U_i,V_i,V'_i,\phi_i)$ (as in Section \ref{sec:degenerate}).  
Note that the quadruple 
corresponding to the initial position is $(Q, \frn, \frn', \phi)$.

We claim that if a root game can be won, then {\em every}
such quadruple encountered over the course of the game is good.
In particular the initial position is good, which, by Lemma
\ref{b-lem:linalg} implies that $\branchproblem \neq 0$.

First, we note that the quadruples $(U_i,V_i,V'_i,\phi_i)$ associated
to a winning position are good.  Indeed, if $\tokens$ is injective, 
then $\phi_i|_{U_i}: U_i \to V'_i$ is an injective linear map,
thus $(U_i,V_i,V'_i,\phi_i)$ is good.

To establish the claim we must show two things:
\begin{enumerate}
\item[(i)] Suppose $(\regions, \tokens)$ is the position of a root game
before a move $[\beta, R_j]$, and $(\regions, \tokens')$ is the position 
after the move.
If all quadruples associated to $(\regions, \tokens')$ are good,
then all quadruples associated to $(\regions, \tokens)$ are good.

\item[(ii)] Suppose $(\regions, \tokens)$ is the position of a root game
before splitting along a splitting subset $A$, and $(\regions', \tokens)$ 
is the position after the splitting.
If all quadruples associated to $(\regions', \tokens)$ are good,
then all quadruples associated to $(\regions, \tokens)$ are good.
\end{enumerate}

\noindent
{\em Proof of (i):}
All quadruples $(U_i,V_i,V'_i,\phi_i)$, $i \neq j$, are unchanged by the 
move $[\beta, R_j]$.  The position $(U_j, V_j, V'_j, \phi_j)$, however,
is changed to $(U_j^1, V_j, V'_j, \phi_j)$, where 
$\Gamma(U_j^1) = \Gamma(U_j) \shift{R_j} \beta$.  By Lemma
\ref{lem:rootsmove}, 
$U_j^1 = \lim_{t \to \infty} \theta_\beta(t) \cdot U_j$,
and thus (i) follows by Proposition \ref{prop:geommove}. \\
\\
\noindent
{\em Proof of (ii):}
Let $S= \langle e_\alpha\ |\ \alpha \in A\rangle$ be the ideal of
$\frn$ corresponding to $A$.  Let $S_i$ be the corresponding submodule
of $V_i$: $S_i = S \cap V_i^0/V_i^1$.  Put $S'_i = \phi_i(S_i)$.
We let $q_i : V_i \to V_i/S_i$,
and $q'_i : V'_i \to V'_i/S_i$,
 denote the quotient maps.

The result of splitting the region $R_i \in \regions$
along $A$ is two regions: $R_i \cap A$ and $R_i \cap A^c$.
Let $(U_i, V_i, V'_i, \phi_i)$ be the quadruple associate to $R_i$.
Then the quadruple associated to $R_i \cap A$
is $(U_i \cap S_i, S_i, S'_i, \phi_i|S_i)$,
 and the quadruple associated $R_i \cap A^c$
$(q(U_i), V_i/S_i, V'_i/S'_i, q'_i \circ \phi_i \circ q_i^{-1})$.  
(The latter, is because
$A$ is a splitting subset (not merely an ideal subset), thus 
$\phi$ respects not just the weight spaces of
$S_i, S'_i$, but also the complementary weight spaces.)
Using Proposition \ref{prop:geomsplit}, (ii) follows.

\end{proof}

\end{document}